\documentclass[AMA,Times1COL]{WileyNJDv5} 

\articletype{Original Article}%

\journal{Preprint, accepted in IJNME}
\startpage{1}

\usepackage{draftwatermark}
\SetWatermarkAngle{90}
\SetWatermarkLightness{.7}
\SetWatermarkHorCenter{1.145\textwidth}
\SetWatermarkFontSize{0.75cm}
\SetWatermarkScale{1}
\SetWatermarkText{{ACCEPTED IN IJNME \;\;\;\;\;\;\;\;\;\; https://doi.org/10.1002/nme.70229}}

\usepackage{amsmath}					
\usepackage{amssymb}					

\mathcode`,="013B						
\mathcode`.="613A						

\usepackage{tikz}
\usepackage{pgfplots}
\usepackage{float}
\usepackage{graphicx,psfrag}

\usepackage{ccaption}
\hangcaption
\captionstyle{\raggedright}

\usepackage{booktabs}
\usepackage{multirow}
\usepackage{tabularray}					

\usepackage{bbold}
\usepackage{accents}
\usepackage{xcolor}

\begin{document}
\title{Development of a Model Order Reduced Arbitrary Lagrangian Eulerian (MORALE) formulation for structures subjected to dynamic moving loads}
\author{Atul Anantheswar$^1$}
\author{Jannick Kehls$^2$}
\author{Ines Wollny$^1$}
\author{Tim Brepols$^2$}
\author{Stefanie Reese$^{2,3}$}
\author{Michael Kaliske$^{1,}$*}

\authormark{\textsc{Anantheswar et al.}}
\address{\orgdiv{$^1$Institute for Structural Analysis}, \orgname{Technische Universit\"at Dresden}, \orgaddress{01062 \state{Dresden}, \country{Germany}}}
\address{\orgdiv{$^2$Institute of Applied Mechanics}, \orgname{RWTH Aachen University}, \orgaddress{52074 \state{Aachen}, \country{Germany}}}
\address{\orgdiv{$^3$University of Siegen}, \orgaddress{57076 \state{Siegen}, \country{Germany}}}
\corres{*Corresponding author, \email{michael.kaliske@tu-dresden.de}}
\presentaddress{Technische Universit\"at Dresden, 01062 Dresden, Germany}


\abstract[Abstract]{Fast and efficient simulation strategies are a basic requirement of technologies such as digital twins. Particularly for roadway infrastructure, recent developments have demonstrated that the Arbitrary Lagrangian Eulerian (ALE) formulation can be utilized to improve computational efficiency, when simulating the response of such pavement structures subjected to moving loads. It is also well established in literature, that Model Order Reduction (MOR) techniques significantly enhance calculation speed. This contribution details the combination of both these tools into a novel Model Order Reduced Arbitrary Lagrangian Eulerian (MORALE) formulation. Both hyperelastic and viscoelastic material models are considered in this work. Transient simulations of pavement structures subjected to moving loads are then carried out with the developed framework, and these show a significant enhancement in computational speed and efficiency over conventional simulation techniques. Also in this work, a comprehensive breakdown of the computational costs involved in using the MORALE formulation is provided and analysed. Such an efficient and fast simulation framework is of vital importance in technologies such as digital twins of roadway infrastructure (like pavements), as it enables engineers to quickly run what-if analyses and make informed decisions about the management of the structure under consideration.}


\maketitle

\section{Introduction}\label{sec:intro}
The contribution highlighted in this manuscript is part of ongoing research in the collaborative research centre SFB/TRR 339. The aim of this project is to establish a digital twin of the road \cite{kaliske2022weg, kaliske2023digitaler}. Such a digital twin of the road necessitates the employment of fast and efficient simulation strategies. One approach enabling fast simulations that can predict the response of pavement structures is the Arbitrary Lagrangian Eulerian (ALE) formulation. \\
ALE approaches are most commonly used in the field of computational fluid dynamics \cite{venkatasubban1995new, souli2000ale, codina2009fixed, basting2017extended}. They are also applied as mesh refinement techniques to counter issues that arise due to extreme mesh distortions \cite{benson1989efficient, bayoumi2004complete, donea2004arbitrary, nazem2009arbitrary, zreid2021ale, berger2022arbitrary}. ALE approaches also find application in moving mesh techniques to describe crack propagation in the field of fracture mechanics \cite{bruno2013fracture,ammendolea2023fatigue,ammendolea2025efficient,ammendolea2025finite}. In recent years, the ALE approach has been used to model the steady state rotation of tires using a fixed finite element mesh \cite{nackenhorst2004ale, ziefle2008numerical}. Subsequent developments have enabled the simulation of the response of pavement structures to steady rolling (constant velocity) wheel loads \cite{wollny2013numerical, wollny2016numerical}. Further research in this direction has allowed the simulation of the transient response of pavements subjected to dynamic wheel loads \cite{anantheswar2023dynamic, anantheswar2023transient}, along with the consideration of inelastic material behaviour \cite{anantheswar2025treatment}. \\
Moreover, another technique to lower the computation time of the simulations is model order reduction (MOR). For finite element simulations, projection-based MOR has been widely used \cite{hesthaven_certified_2016, benner_survey_2015}. The main idea of projection-based MOR is to reduce the size of the system of equations by projecting the high-dimensional system onto a low-dimensional subspace. This is achieved by using a set of basis functions that span the low-dimensional subspace. The reduced system of equations can then be solved more efficiently than the original system. The reduced order model (ROM) can be constructed using various techniques, such as proper orthogonal decomposition (POD) \cite{benner_survey_2015}, balanced truncation \cite{gugercin_survey_2004}, and modal basis method \cite{spiess_reduction_2005}. In this work, we focus on the POD method for constructing the ROM. It has been shown that the POD delivers accurate approximations for nonlinear problems in solid mechanics \cite{radermacher_comparison_2013} and has been extended for various applications such as component-based MOR \cite{ritzert2025component}, adaptively sub-structured elastoplastic simulations \cite{radermacher2014model}, and reduced-order modelling of gradient-extended damage \cite{zhang2025multi}. Furthermore, the POD serves as a basis for most hyper-reduction techniques such as the discrete empirical interpolation method \cite{chaturantabut2010nonlinear, radermacher2016pod, kehls2023reduced} or energy-conserving sampling and weighting \cite{an2008optimizing, farhat2014dimensional}. The POD is a data-driven approach that uses snapshots of the solution to construct a low-dimensional basis that captures the dominant features of the solution.\\
Recent developments showcase the coupling of ALE and MOR in the simulation of rolling tires in contact with the road \cite{de2019application, de2019development, wang2024nonsmooth}. In these works, the behaviour of the tire is the main focus, and the combination of ALE and MOR techniques has been demonstrated to improve computational efficiency. However, in the present work, the combination of ALE and MOR frameworks to simulate the transient response of longitudinally homogeneous pavement structures subjected to moving loads, is detailed. \\
This manuscript is structured as follows. In Section \ref{sec:kinematics}, the reasoning behind the adoption of the ALE approach to simulate the response of pavement structures is provided with a kinematics background. Then, the balance laws which are the basis for an implementation into a finite element framework are provided in Section \ref{sec:ALE_FE}, specifically keeping the ALE reference configuration in mind. Details on the advection of the solution field as well as material point quantities are given in Section \ref{sec:ALE_updates}. Next, Section \ref{sec:MOR} describes the application of model order reduction techniques to the assembled system of equations which includes contributions from the ALE implementation. Further, numerical studies are carried out and the results are presented in Section \ref{sec:numex}, highlighting the benefits of the combined MORALE formulation. Finally, the conclusions drawn from the studies are summarized in Section \ref{sec:conc}.

\section{Kinematics in the ALE formulation}\label{sec:kinematics}
\begin{figure}[h!]
	\centering
	\def\svgwidth{\textwidth}
	\begingroup%
	\makeatletter%
	\providecommand\color[2][]{%
		\errmessage{(Inkscape) Color is used for the text in Inkscape, but the package 'color.sty' is not loaded}%
		\renewcommand\color[2][]{}%
	}%
	\providecommand\transparent[1]{%
		\errmessage{(Inkscape) Transparency is used (non-zero) for the text in Inkscape, but the package 'transparent.sty' is not loaded}%
		\renewcommand\transparent[1]{}%
	}%
	\providecommand\rotatebox[2]{#2}%
	\newcommand*\fsize{\dimexpr\f@size pt\relax}%
	\newcommand*\lineheight[1]{\fontsize{\fsize}{#1\fsize}\selectfont}%
	\ifx\svgwidth\undefined%
	\setlength{\unitlength}{1723.01075372bp}%
	\ifx\svgscale\undefined%
	\relax%
	\else%
	\setlength{\unitlength}{\unitlength * \real{\svgscale}}%
	\fi%
	\else%
	\setlength{\unitlength}{\svgwidth}%
	\fi%
	\global\let\svgwidth\undefined%
	\global\let\svgscale\undefined%
	\makeatother%
	\begin{picture}(1,0.59384814)%
		\lineheight{1}%
		\setlength\tabcolsep{0pt}%
		\put(0,0){\includegraphics[width=\unitlength,page=1]{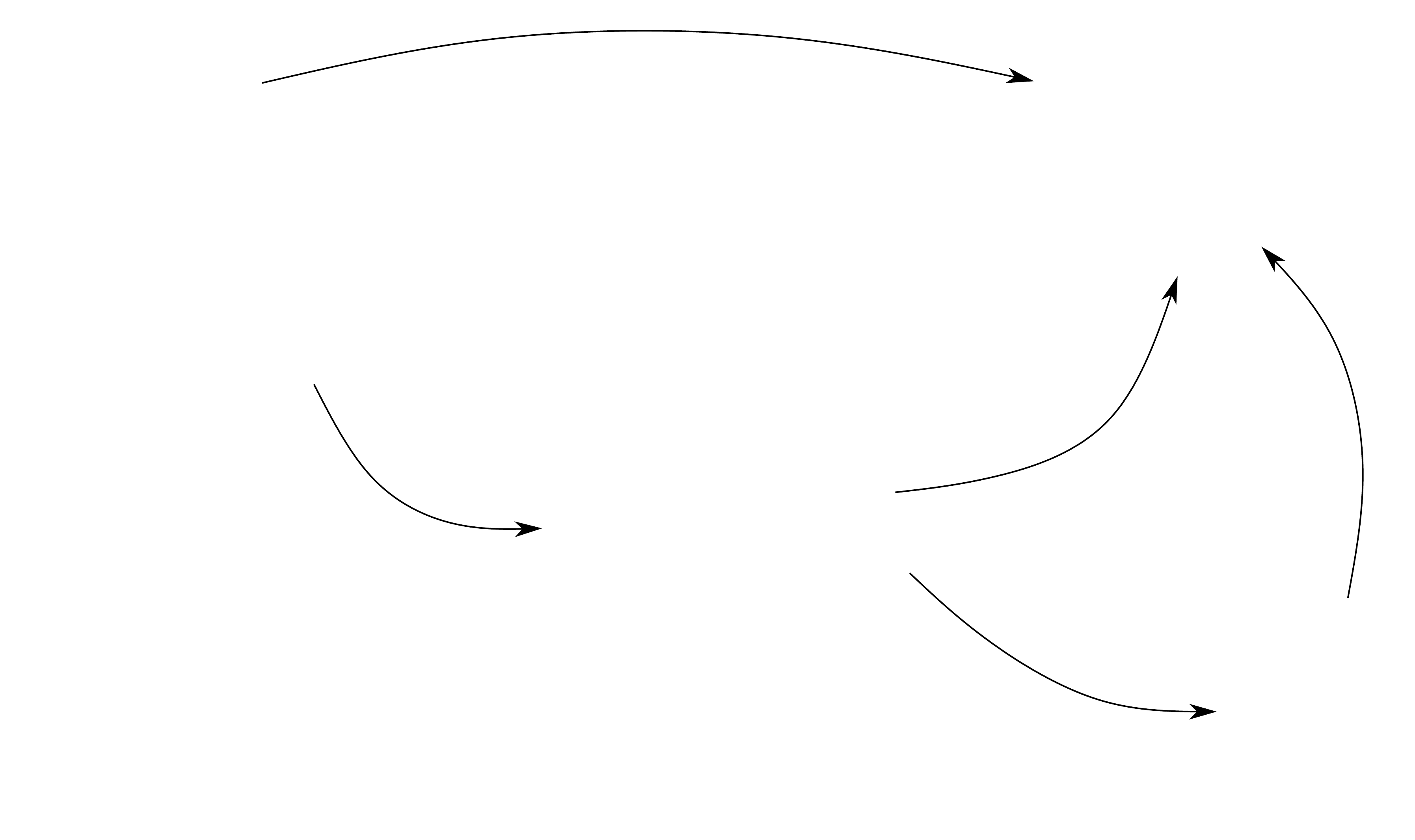}}%
		\put(0.06190534,0.58500308){\makebox(0,0)[lt]{\lineheight{1.25}\smash{\begin{tabular}[t]{l}initial configuration\end{tabular}}}}%
		\put(0.50327362,0.09707869){\makebox(0,0)[lt]{\lineheight{1.25}\smash{\begin{tabular}[t]{l}ALE reference \\\end{tabular}}}}%
		\put(0.73099383,0.58551386){\makebox(0,0)[lt]{\lineheight{1.25}\smash{\begin{tabular}[t]{l}current configuration\end{tabular}}}}%
		\put(0.80048304,0.02830044){\makebox(0,0)[lt]{\lineheight{1.25}\smash{\begin{tabular}[t]{l}inelastic intermediate\end{tabular}}}}%
		\put(0,0){\includegraphics[width=\unitlength,page=2]{01_ALE_Kinematics_svg-tex.pdf}}%
		\put(0.45121656,0.59032053){\makebox(0,0)[lt]{\lineheight{1.25}\smash{\begin{tabular}[t]{l}$\phi$\end{tabular}}}}%
		\put(0,0){\includegraphics[width=\unitlength,page=3]{01_ALE_Kinematics_svg-tex.pdf}}%
		\put(0.29403234,0.25326005){\makebox(0,0)[lt]{\lineheight{1.25}\smash{\begin{tabular}[t]{l}$\chi$\end{tabular}}}}%
		\put(0.72765414,0.28829828){\makebox(0,0)[lt]{\lineheight{1.25}\smash{\begin{tabular}[t]{l}$\hat{\phi}$\end{tabular}}}}%
		\put(0.45018146,0.53922771){\makebox(0,0)[lt]{\lineheight{1.25}\smash{\begin{tabular}[t]{l}$\boldsymbol{\underline{F}}$\end{tabular}}}}%
		\put(0.25913249,0.20450958){\makebox(0,0)[lt]{\lineheight{1.25}\smash{\begin{tabular}[t]{l}$\boldsymbol{\underline{R}}$\end{tabular}}}}%
		\put(0.77413059,0.25131499){\makebox(0,0)[lt]{\lineheight{1.25}\smash{\begin{tabular}[t]{l}$\boldsymbol{\underline{\hat{F}}}$\end{tabular}}}}%
		\put(0.72980409,0.14521912){\makebox(0,0)[lt]{\lineheight{1.25}\smash{\begin{tabular}[t]{l}$\boldsymbol{\underline{\hat{F}}}^{\text{i}}$\end{tabular}}}}%
		\put(0.90033343,0.30314568){\makebox(0,0)[lt]{\lineheight{1.25}\smash{\begin{tabular}[t]{l}$\boldsymbol{\underline{\hat{F}}}^{\text{e}}$\end{tabular}}}}%
		\put(0,0){\includegraphics[width=\unitlength,page=4]{01_ALE_Kinematics_svg-tex.pdf}}%
		\put(0.05524438,0.39887465){\makebox(0,0)[lt]{\lineheight{1.25}\smash{\begin{tabular}[t]{l}$\boldsymbol{X}$\end{tabular}}}}%
		\put(0.29517418,0.48691171){\makebox(0,0)[lt]{\lineheight{1.25}\smash{\begin{tabular}[t]{l}$\boldsymbol{u}$\end{tabular}}}}%
		\put(0.3347482,0.39714822){\makebox(0,0)[lt]{\lineheight{1.25}\smash{\begin{tabular}[t]{l}$\boldsymbol{x}$\end{tabular}}}}%
		\put(0.33783149,0.31383369){\makebox(0,0)[lt]{\lineheight{1.25}\smash{\begin{tabular}[t]{l}$\boldsymbol{u}_{\text{rig}}$\end{tabular}}}}%
		\put(0.63763945,0.33260453){\makebox(0,0)[lt]{\lineheight{1.25}\smash{\begin{tabular}[t]{l}$\boldsymbol{\varphi}$\end{tabular}}}}%
		\put(0.39853783,0.16359431){\makebox(0,0)[lt]{\lineheight{1.25}\smash{\begin{tabular}[t]{l}$\boldsymbol{\chi}$\end{tabular}}}}%
		\put(0.57550711,0.35596208){\makebox(0,0)[lt]{\lineheight{1.25}\smash{\begin{tabular}[t]{l}$\boldsymbol{\hat{u}}$\end{tabular}}}}%
		\put(0.14202261,0.13471803){\makebox(0,0)[lt]{\lineheight{1.25}\smash{\begin{tabular}[t]{l}$\boldsymbol{\chi}_{\text{0}}^{\text{ALE}}$\end{tabular}}}}%
		\put(-0.00038542,0.20280504){\makebox(0,0)[lt]{\lineheight{1.25}\smash{\begin{tabular}[t]{l}$\boldsymbol{e}_{\text{i}}^{\text{Lagrange}}$\end{tabular}}}}%
		\put(0.38002239,0.01796471){\makebox(0,0)[lt]{\lineheight{1.25}\smash{\begin{tabular}[t]{l}$\boldsymbol{e}_{\text{i}}^{\text{ALE}}$\end{tabular}}}}%
		\put(0,0){\includegraphics[width=\unitlength,page=5]{01_ALE_Kinematics_svg-tex.pdf}}%
		\put(0.51239129,0.06821123){\makebox(0,0)[lt]{\lineheight{1.25}\smash{\begin{tabular}[t]{l}configuration\end{tabular}}}}%
		\put(0.84884254,0.00096581){\makebox(0,0)[lt]{\lineheight{1.25}\smash{\begin{tabular}[t]{l}configuration\end{tabular}}}}%
	\end{picture}%
	\endgroup%
	\caption{ALE kinematics \cite{anantheswar2025treatment}.}
	\label{fig:ALE_kinematics}
\end{figure}

Consider the deformation of a body from the initial configuration to the current configuration as shown in Figure \ref{fig:ALE_kinematics}. Conventional methods commonly assume that the reference configuration is identical to the initial configuration and both share the same stationary reference frame $\boldsymbol{e}_{\text{i}}^{\text{Lagrange}}$. However, for the ALE approach, a reference frame that moves with time $\boldsymbol{e}_{\text{i}}^{\text{ALE}}$ is introduced. For the particular case of structures subjected to moving loads, the ALE reference frame can be moved with the same velocity as that of the moving load. This moving reference frame then offers a unique perspective, where the load appears stationary but the material of the structure (pavements, gantry crane girders etc.) appears to flow beneath the load \cite{wollny2013numerical}. This comes with the advantage that only the relevant region (domain of interest) around the load needs to be analysed. Hence, the computational costs associated with discretizing and simulating regions far away from the load are avoided. When it comes to structures like pavements subjected to moving vehicle loads, this irrelevant domain can be significantly large, as the entire length of the structure in the path of the moving load would need to be analysed. It is also worth mentioning that the deformation from the ALE reference configuration can be further decomposed to accommodate an inelastic intermediate configuration, thereby enabling the simulation of inelastic structures \cite{anantheswar2025treatment}.

From Figure \ref{fig:ALE_kinematics}, the current position $\boldsymbol{x}$ of a material point can be written in terms of the position vectors ($\boldsymbol{X}$, $\boldsymbol{\chi}$ and $\boldsymbol{\chi}_{\text{0}}^{\text{ALE}}$) and displacements ($\boldsymbol{u}$, $\boldsymbol{u}_{\text{rig}}$ and $\hat{\boldsymbol{u}}$) in various reference frames as
\begin{equation}\label{eq:current_position}
	\boldsymbol{x} = \boldsymbol{X} + \boldsymbol{u} = \boldsymbol{X} + \boldsymbol{u}_{\textnormal{rig}} + \boldsymbol{\hat{u}} = \boldsymbol{\chi}_{\textnormal{0}}^{\textnormal{ALE}} + \boldsymbol{\chi} + \boldsymbol{\hat{u}} \;\;.
\end{equation}
It is to be noted that for the pavement structures under consideration in this manuscript, the rigid body displacements $\boldsymbol{u}_{\text{rig}} = \boldsymbol{0}$, because there are no rigid translations or rotations. Therefore, at time $t = 0$, the ALE reference frame and the initial reference frame coincide. However, the ALE reference frame moves with time along with the moving load. Hence, discretization is carried out with respect to positions $\boldsymbol{\chi}$ in the moving ALE reference frame, and a new portion of the overall structure enters the mesh at each point in time. Further, the velocity with which the material of the structure appears to flow through the mesh is termed `guiding velocity', denoted $\boldsymbol{w}$, and expressed as
\begin{equation}\label{eq:wpav_def}
	\boldsymbol{w} = \left.{\frac{\partial\boldsymbol{\chi}}{\partial t}} \right|_{\boldsymbol{X}} = \left.{\frac{\partial\left(\boldsymbol{X} + \boldsymbol{u}_{\textnormal{rig}} - \boldsymbol{\chi}_{\textnormal{0}}^{\textnormal{ALE}}\right)}{\partial t}} \right|_{\boldsymbol{X}} = \left.-{\frac{\partial\boldsymbol{\chi}_{\textnormal{0}}^{\textnormal{ALE}}}{\partial t}} \right|_{\boldsymbol{X}} \;\;.
\end{equation}

\section{Dynamic ALE formulation in a Finite Element framework}\label{sec:ALE_FE}
The implementation of the ALE formulation into a finite element framework \cite{taylor_feap_2020} is as per previous research \cite{anantheswar2023dynamic,anantheswar2025treatment}, and is briefly summarized here. The basis for a finite element implementation of a mechanical problem is the balance of linear momentum. For the transient case in the ALE reference configuration, the strong form of this balance law can be expressed \cite{nackenhorst2004ale,wollny2013numerical} as
\begin{equation}\label{eq:lin_mom}
	\hat{\rho}\boldsymbol{\dot{v}} = \textnormal{Div}\boldsymbol{\hat{\underline{P}}} + \hat{\rho}\boldsymbol{b}  \;\;,
\end{equation}
where $\hat{\rho}$ refers to material density, $\boldsymbol{\hat{\underline{P}}}$ denotes the first Piola-Kirchhoff stress, and $\boldsymbol{b}$ refers to the applied body forces per unit mass. Throughout this manuscript, quantities with a hat $\hat{\square}$ denote that they are in the ALE reference configuration. Further, the weak form is obtained by multiplying a test function $\boldsymbol{\eta}$, and integrating Equation (\ref{eq:lin_mom}) over the volume of the body in the ALE reference configuration with the help of partial integration and the Gauss' theorem. The weak form is
\begin{equation}\label{eq:lin_mom_weak}
	\int_{\chi(B)}\hat{\rho}\boldsymbol{\dot{v}}\cdot\boldsymbol{\eta}d\hat{V} + \int_{\chi(B)}\boldsymbol{\hat{\underline{P}}}:\textnormal{Grad}\boldsymbol{\eta}d\hat{V} = \int_{\chi(B)}\hat{\rho}\boldsymbol{b}\cdot\boldsymbol{\eta}d\hat{V} + \int_{\partial\chi(B)}\boldsymbol{\hat{T}}\cdot\boldsymbol{\eta}d\hat{A}  \;\;.
\end{equation}
Here, $\boldsymbol{\hat{T}}$ refers to the externally applied surface traction, and the gradient operator $\text{Grad}\square = \frac{\partial \square}{\partial \boldsymbol{\chi}}$ is the gradient in the ALE reference configuration. Additionally, $\chi(B)$ refers to the domain of the body and $\partial\chi(B)$ denotes its surface, both in the ALE reference frame. For the dynamic ALE formulation, the first term (representing inertial contributions) in Equation (\ref{eq:lin_mom_weak}) is different when compared to conventional formulations. This is because in the ALE reference frame, the material appears to flow through the mesh, and the advection of velocity and acceleration terms (which are tied to the material points) needs to be considered. These are written as
\begin{equation}\label{eq:vel_con}
	\boldsymbol{v} (\boldsymbol{x},t) = \left.\frac{\partial\boldsymbol{x}}{\partial t} \right|_{\boldsymbol{X}} = \left.\frac{\partial\boldsymbol{x}}{\partial t} \right|_{\boldsymbol{\chi}} + \textnormal{Grad}\boldsymbol{x}\cdot\boldsymbol{w} \;\;,
\end{equation}
and
\begin{equation}\label{eq:acc_con}
	\boldsymbol{\dot{v}} (\boldsymbol{x},t) = \left.\frac{\partial\boldsymbol{v} (\boldsymbol{x},t)}{\partial t} \right|_{\boldsymbol{X}} = \left.\frac{\partial\boldsymbol{v} (\boldsymbol{x},t)}{\partial t} \right|_{\boldsymbol{\chi}} + \textnormal{Grad}\boldsymbol{v} (\boldsymbol{x},t)\cdot\boldsymbol{w} \;\;,
\end{equation}
for the velocity and acceleration, respectively \cite{donea2004arbitrary,nackenhorst2004ale,wollny2013numerical,anantheswar2023dynamic}. The position vector in the current configuration is $\boldsymbol{x} = \boldsymbol{\chi}_\textnormal{0}^\textnormal{ALE} + \boldsymbol{\varphi}$, and so, Equation (\ref{eq:vel_con}) can be written in terms of $\boldsymbol{\varphi}$ as
\begin{equation}\label{eq:v_interms_phi}
	\boldsymbol{v} (\boldsymbol{x},t) = \left.\frac{\partial\boldsymbol{\chi}_{\textnormal{0}}^{\textnormal{ALE}}}{\partial t} \right|_{\boldsymbol{\chi}} + \left.\frac{\partial\boldsymbol{\varphi}}{\partial t} \right|_{\boldsymbol{\chi}} + \textnormal{Grad}\boldsymbol{\varphi}\cdot\boldsymbol{w} \;\;,
\end{equation}
because $\textnormal{Grad}\boldsymbol{\chi}_{\textnormal{0}}^{\textnormal{ALE}} = \boldsymbol{\underline{0}}$. Further, substituting Equation (\ref{eq:v_interms_phi}) into the Equation (\ref{eq:acc_con}), it is possible to write
\begin{equation}\label{eq:a_interms_phi}
	\begin{aligned}
		\boldsymbol{\dot{v}} (\boldsymbol{x},t) = &
		\left.\frac{\partial^{2}\boldsymbol{\chi}_{\textnormal{0}}^{\textnormal{ALE}}}{\partial t^{2}}\right|_{\boldsymbol{\chi}} +
		\left.\frac{\partial^{2}\boldsymbol{\varphi}}{\partial t^{2}}\right|_{\boldsymbol{\chi}} + \\
		&2\cdot\textnormal{Grad}\left(\left.\frac{\partial\boldsymbol{\varphi}}{\partial t}\right|_{\boldsymbol{\chi}}\right)\cdot\boldsymbol{w} +
		\textnormal{Grad}\boldsymbol{\varphi}\cdot\left.\frac{\partial\boldsymbol{w}}{\partial t}\right|_{\boldsymbol{\chi}} + \\
		&\textnormal{Grad}\left(\textnormal{Grad}\boldsymbol{\varphi}\cdot\boldsymbol{w}\right)\cdot\boldsymbol{w} \;\;.
	\end{aligned}
\end{equation}
Finally, the above expression for the acceleration is substituted into the inertia term of Equation (\ref{eq:lin_mom_weak}), and the following expression for the inertia term can be obtained \cite{anantheswar2023dynamic}
\begin{equation}\label{eq:ma_term_ale_pav}
	\begin{aligned}
		\int_{\chi(B)}\hat{\rho}\boldsymbol{\dot{v}}\cdot\boldsymbol{\eta}d\hat{V} = &
		\int_{\chi(B)}\hat{\rho}\left.\frac{\partial^{2}\boldsymbol{\chi}_{\textnormal{0}}^{\textnormal{ALE}}}{\partial t^{2}}\right|_{\boldsymbol{\chi}}\cdot\boldsymbol{\eta}d\hat{V} +
		\int_{\chi(B)}\hat{\rho}\left.\frac{\partial^{2}\boldsymbol{\varphi}}{\partial t^{2}}\right|_{\boldsymbol{\chi}}\cdot\boldsymbol{\eta}d\hat{V} + \\
		&2\cdot\int_{\chi(B)}\hat{\rho}\;\textnormal{Grad}\left(\left.\frac{\partial\boldsymbol{\varphi}}{\partial t}\right|_{\boldsymbol{\chi}}\right)\cdot\boldsymbol{w}\cdot\boldsymbol{\eta}d\hat{V} + \\
		&\int_{\chi(B)}\hat{\rho}\;\textnormal{Grad}\boldsymbol{\varphi}\cdot\left.\frac{\partial\boldsymbol{w}}{\partial t}\right|_{\boldsymbol{\chi}}\cdot\boldsymbol{\eta}d\hat{V} + \\
		&\int_{\chi(B)}\hat{\rho}\;\left(\textnormal{Grad}\left(\textnormal{Grad}\boldsymbol{\varphi}\right)\cdot\boldsymbol{w}\right)\cdot\boldsymbol{w}\cdot\boldsymbol{\eta}d\hat{V} \;\;.
	\end{aligned}
\end{equation}
Note that if the guiding velocity $\boldsymbol{w}$ varies in the mesh, additional terms would arise in the above Equation (\ref{eq:ma_term_ale_pav}), and would also need to be considered \cite{nackenhorst2004ale,anantheswar2023dynamic}. The linearization and discretization of Equation (\ref{eq:ma_term_ale_pav}) \cite{anantheswar2023dynamic, anantheswar2025treatment} is summarized in Table \ref{table:summary_disc2}.
\begin{table}[h!]
	\caption{Summary of the linearized terms obtained from Equation (\ref{eq:ma_term_ale_pav}).}
	\begin{center}
		\begin{tblr}{colspec = {X[c,m]X[l,m]X[c,m]X[l,m]},
				cell{1}{1} = {c=1}{c}
			}
			\hline
			\hline
			\textbf{No.} & \textbf{tangent} & \textbf{increment} & \textbf{internal force} \\
			\hline
			1 & $\boldsymbol{\underline{0}}$ & - & $-\int_{\chi^{e}}\hat{\rho}\boldsymbol{\underline{H}}^{T}\cdot\left.\frac{\partial\boldsymbol{w}}{\partial t}\right|_{\chi}d\hat{V}^{e}$ \\
			2 & $\int_{\chi^{e}}\hat{\rho}\boldsymbol{\underline{H}}^{T}\cdot\boldsymbol{\underline{H}} d\hat{V}^{e}$ & $\Delta\frac{\widetilde{\partial^{2}\boldsymbol{\varphi}}}{\partial t^{2}}$ & $\int_{\chi^{e}}\hat{\rho}\boldsymbol{\underline{H}}^{T}\cdot\boldsymbol{\underline{H}}\cdot\left.\frac{\widetilde{\partial^{2}\boldsymbol{\varphi}}}{\partial t^{2}}\right|_{\chi}d\hat{V}^{e}$ \\
			3 & $2\int_{\chi^{e}}\hat{\rho}\boldsymbol{\underline{H}}^{T}\cdot\boldsymbol{\underline{A}} d\hat{V}^{e}$ & $\Delta\frac{\widetilde{\partial\boldsymbol{\varphi}}}{\partial t}$ & $2\int_{\chi^{e}}\hat{\rho}\boldsymbol{\underline{H}}^{T}\cdot\boldsymbol{\underline{A}}\cdot\left.\frac{\widetilde{\partial\boldsymbol{\varphi}}}{\partial t}\right|_{\chi}d\hat{V}^{e}$ \\
			4 & $\int_{\chi^{e}}\hat{\rho}\boldsymbol{\underline{H}}^{T}\cdot\boldsymbol{\underline{A'}}d\hat{V}^{e}$ & $\Delta\widetilde{\boldsymbol{\varphi}}$ & $\int_{\chi^{e}}\hat{\rho}\boldsymbol{\underline{H}}^{T}\cdot\boldsymbol{\underline{A'}}\cdot\widetilde{\boldsymbol{\varphi}}d\hat{V}^{e}$ \\
			5 & $\int_{\chi^{e}}\hat{\rho}\boldsymbol{\underline{H}}^{T}\cdot\boldsymbol{\underline{A}}''d\hat{V}^{e}$ & $\Delta\widetilde{\boldsymbol{\varphi}}$ & $\int_{\chi^{e}}\hat{\rho}\boldsymbol{\underline{H}}^{T}\cdot\boldsymbol{\underline{A}}''\cdot\widetilde{\boldsymbol{\varphi}}d\hat{V}^{e}$ \\
			\hline
			\hline
		\end{tblr}
	\end{center}
	\label{table:summary_disc2}
\end{table}

In Table \ref{table:summary_disc2}, the term $\boldsymbol{\underline{H}}$ is an orderly arrangement of the shape functions, given by
\begin{equation}\label{eq:Hmatrix}
	\boldsymbol{\underline{H}}(\boldsymbol{\xi}) = \left[\begin{array}{ccccccc}
		N_{1}(\boldsymbol{\xi}) & 0 & 0 & \cdots & N_{k}(\boldsymbol{\xi}) & 0 & 0 \\
		0 & N_{1}(\boldsymbol{\xi}) & 0 & \cdots & 0 & N_{k}(\boldsymbol{\xi}) & 0 \\
		0 & 0 & N_{1}(\boldsymbol{\xi}) & \cdots & 0 & 0 & N_{k}(\boldsymbol{\xi})
	\end{array}\right] \;\;,
\end{equation}
where $N_{k}(\boldsymbol{\xi})$ refers to the shape function associated with iso-parametric coordinates $\boldsymbol{\xi}$ and node $k$ of a given finite element. The matrix $\boldsymbol{\underline{H}}$ can be utilized for interpolation from nodal values of any quantity $\widetilde{\square}$ to the value at integration points within a finite element. Also, in Table \ref{table:summary_disc2}, $\boldsymbol{\underline{A}}$ denotes the matrix
\begin{equation}\label{eq:Amatrix}
	\boldsymbol{\underline{A}} = \left[\begin{array}{ccccccc}
		N_{1,i}w_{i} & 0 & 0 & \cdots & N_{k,i}w_{i} & 0 & 0 \\
		0 & N_{1,i}w_{i} & 0 & \cdots & 0 & N_{k,i}w_{i} & 0 \\
		0 & 0 & N_{1,i}w_{i} & \cdots & 0 & 0 & N_{k,i}w_{i}
	\end{array}\right] \;\;,
\end{equation}
where $N_{k,i}w_{i}$ is summed over index $i$. $N_{k,i}$ is the gradient of the shape function associated with node $k$ along $i$-direction, and $w_{i}$ is the component of the guiding velocity along $i$-direction. The term $\boldsymbol{\underline{A}}'$ in Table \ref{table:summary_disc2} refers to the matrix
\begin{equation}\label{eq:A'matrix}
	\boldsymbol{\underline{A}}' = \left[\begin{array}{ccccccc}
		N_{1,i}{\dot{w}}_{i} & 0 & 0 & \cdots & N_{k,i}{\dot{w}}_{i} & 0 & 0 \\
		0 & N_{1,i}{\dot{w}}_{i} & 0 & \cdots & 0 & N_{k,i}{\dot{w}}_{i} & 0 \\
		0 & 0 & N_{1,i}{\dot{w}}_{i} & \cdots & 0 & 0 & N_{k,i}{\dot{w}}_{i}
	\end{array}\right] \;\;,
\end{equation}
where ${\dot{w}}_{i}$ is the component of $\left.\frac{\partial\boldsymbol{w}}{\partial t}\right|_{\chi}$ along the $i$-direction. Since all nodes in the mesh have the same guiding velocity at any given point of time, the matrix $\boldsymbol{\underline{A'}}$ can be calculated in a simple manner in structures like pavements. However, in structures such as rolling wheels, there is typically a non-zero angular guiding velocity. Therefore, the calculation of $\boldsymbol{\underline{A}'}$ requires the time history data of the guiding velocity. Further, the term $\boldsymbol{\underline{A}}''$ from Table \ref{table:summary_disc2} is a matrix given by
	\begin{equation}\label{eq:A''matrix}
		\begin{aligned}
			\boldsymbol{\underline{A}}'' =
			\left[\begin{array}{ccccccc}
				N_{1,ij}{w}_{j}{w}_{i}&0&0&\small{\cdots}&N_{k,ij}{w}_{j}{w}_{i}&0&0 \\
				0&N_{1,ij}{w}_{j}{w}_{i}&0&\small{\cdots}&0&N_{k,ij}{w}_{j}{w}_{i}&0 \\
				0&0&N_{1,ij}{w}_{j}{w}_{i}&\small{\cdots}&0&0&N_{k,ij}{w}_{j}{w}_{i}
			\end{array}\right] \;\;,
		\end{aligned}
	\end{equation}
where $N_{i,jk}$ is the second derivative of the shape function corresponding to node $i$ in relation to global coordinates $\boldsymbol{\chi}$.

For the numerical studies in this work, the stress response is obtained using the St. Venant Kirchhoff hyperelastic model \cite{anantheswar2023dynamic}, as well as a nonlinear viscoelastic model based on the Neo-Hookean hyperelastic model \cite{anantheswar2025treatment}. The linearization and discretization of these models is well established in literature \cite{holzapfel2002nonlinear,wriggers2008nonlinear,zopf2015continuum}. The contributions from the inertia term and stress response from the material models can be assembled into the global system of equations in the ALE reference frame, which has to be solved in every iteration $(j)$ of the global Newton-Raphson scheme and can be written as
\begin{equation}
    \label{eq:soe}
    \underline{\boldsymbol{K}}_2^{(j)} \Delta \frac{\partial^2 \boldsymbol{\varphi}}{\partial t^2}+\underline{\boldsymbol{K}}_3^{(j)} \Delta \frac{\partial \boldsymbol{\varphi}}{\partial t} + \left(\underline{\boldsymbol{K}}_4^{(j)}+\underline{\boldsymbol{K}}_5^{(j)} + \underline{\boldsymbol{K}}_T^{(j)} \right) \Delta \boldsymbol{\varphi}=\mathbf{f}_{e x t}^{(j+1)}-\sum_{i=1}^5 \mathbf{f}_i^{(j)} - \mathbf{f}_{int}^{(j)} \;\;.
\end{equation}
In the above Equation (\ref{eq:soe}), $\underline{\boldsymbol{K}}_2$ to $\underline{\boldsymbol{K}}_5$ are the assembled ALE tangent matrices and $\mathbf{f}_1$ to $\mathbf{f}_5$ are the assembled internal forces coming from the ALE formulation as described in Table \ref{table:summary_disc2}. The matrix $\underline{\boldsymbol{K}}_T$ is the global tangential stiffness matrix resulting from the material response, and $\mathbf{f}_{int}$ is the corresponding global internal force vector. Additionally, $\mathbf{f}_{e x t}$ denotes the assembled vector of external forces. Using the Newmark-beta method, Equation (\ref{eq:soe}) can be reformulated as
\begin{equation}\label{eq:KTdyn}
    \underline{\boldsymbol{K}}_{T, dyn} \Delta \boldsymbol{\varphi} = \mathbf{G} \;\;.
\end{equation}
The global dynamic tangential stiffness matrix is computed as
\begin{equation}
    \underline{\boldsymbol{K}}_{T, dyn} = \underline{\boldsymbol{K}}_{T} + \frac{1}{\beta \Delta t^2} \underline{\boldsymbol{K}}_{2} + \frac{\gamma}{\beta \Delta t} \underline{\boldsymbol{K}}_{3} + \underline{\boldsymbol{K}}_{4} + \underline{\boldsymbol{K}}_{5} \;\;,
\end{equation}
where $\beta$ and $\gamma$ are time integration parameters and $\Delta t$ denotes the time increment. Further considerations for the dynamic ALE formulation include the solution field (displacement, velocity and acceleration) update, and history variable update procedures. Theses aspects are also as per previous work \cite{anantheswar2025treatment,am3p_paper_atul}, and for the sake of completeness, are briefly detailed in the next Section \ref{sec:ALE_updates}.

\section{Update procedures required in the dynamic ALE formulation}\label{sec:ALE_updates}
In transient numerical simulations, time integration is necessary. This means that the solution field (displacement, velocity and acceleration) of any time step depends on that at the previous time step (for the Newmark scheme used here). Since the material of the structure appears to flow through the mesh in the ALE reference frame, an update of the solution fields to values of the new material points is necessary. The algorithm used for this process is provided \cite{anantheswar2025treatment} in Table \ref{table:disp_upd}.

\begin{table}[h!]
	\caption{Algorithm used to perform the solution field update between any two consecutive time steps \cite{anantheswar2025treatment}.}
	\begin{center}
		\begin{tabular}{p{0.1\textwidth}p{0.8\textwidth}}
			\hline
			\hline
			\textbf{Start} & \textbf{$t$ = $t_n$; Solution fields are known} \\
			\hline
			1. & Determine coordinates of the new material points that would move into the mesh nodes at $t_{n+1} = t_n + \Delta t$, i.e. $\boldsymbol{\chi}_{n+1} = \boldsymbol{\chi}_{n} - \boldsymbol{w}\Delta t$ \\
			2. & For each mesh node, loop over all elements and identify the element containing the corresponding new material point \\
			3. & Find the isoparametric coordinates ($\boldsymbol{\xi}$) within the identified element corresponding to each new material point \\
			4. & Ensure that these coordinates lie within the identified finite element $\textnormal{abs}[\boldsymbol{\xi}]_{i}<1$; $i$ - direction \\
			5. & Formulate shape functions $\boldsymbol{N}$ to project values from nodes of identified element to the new material point \\
			6. & Interpolate for the value of the solution field at each material point using $\left.\boldsymbol{u}\right|_{t_{n+1}} = \sum_{k}\boldsymbol{N}_{k}\left.\boldsymbol{u}_{k}\right|_{t_n}$; $k$ - mesh nodes \\
			7. & Overwrite $\left.\boldsymbol{u}\right|_{t_n}$ with the interpolated values $\left.\boldsymbol{u}\right|_{t_{n+1}}$ \\
			\hline
			\textbf{End} & \textbf{Solution fields updated for start of next step $t$ = $t_n + \Delta t$} \\
			\hline
			\hline
		\end{tabular}
	\end{center}
	\label{table:disp_upd}
\end{table}

\begin{figure}[h!]
	\centering
	\def\svgwidth{\textwidth}
	\begingroup%
	\makeatletter%
	\providecommand\color[2][]{%
		\errmessage{(Inkscape) Color is used for the text in Inkscape, but the package 'color.sty' is not loaded}%
		\renewcommand\color[2][]{}%
	}%
	\providecommand\transparent[1]{%
		\errmessage{(Inkscape) Transparency is used (non-zero) for the text in Inkscape, but the package 'transparent.sty' is not loaded}%
		\renewcommand\transparent[1]{}%
	}%
	\providecommand\rotatebox[2]{#2}%
	\newcommand*\fsize{\dimexpr\f@size pt\relax}%
	\newcommand*\lineheight[1]{\fontsize{\fsize}{#1\fsize}\selectfont}%
	\ifx\svgwidth\undefined%
	\setlength{\unitlength}{4092.33456109bp}%
	\ifx\svgscale\undefined%
	\relax%
	\else%
	\setlength{\unitlength}{\unitlength * \real{\svgscale}}%
	\fi%
	\else%
	\setlength{\unitlength}{\svgwidth}%
	\fi%
	\global\let\svgwidth\undefined%
	\global\let\svgscale\undefined%
	\makeatother%
	\begin{picture}(1,0.65509609)%
		\lineheight{1}%
		\setlength\tabcolsep{0pt}%
		\put(0,0){\includegraphics[width=\unitlength,page=1]{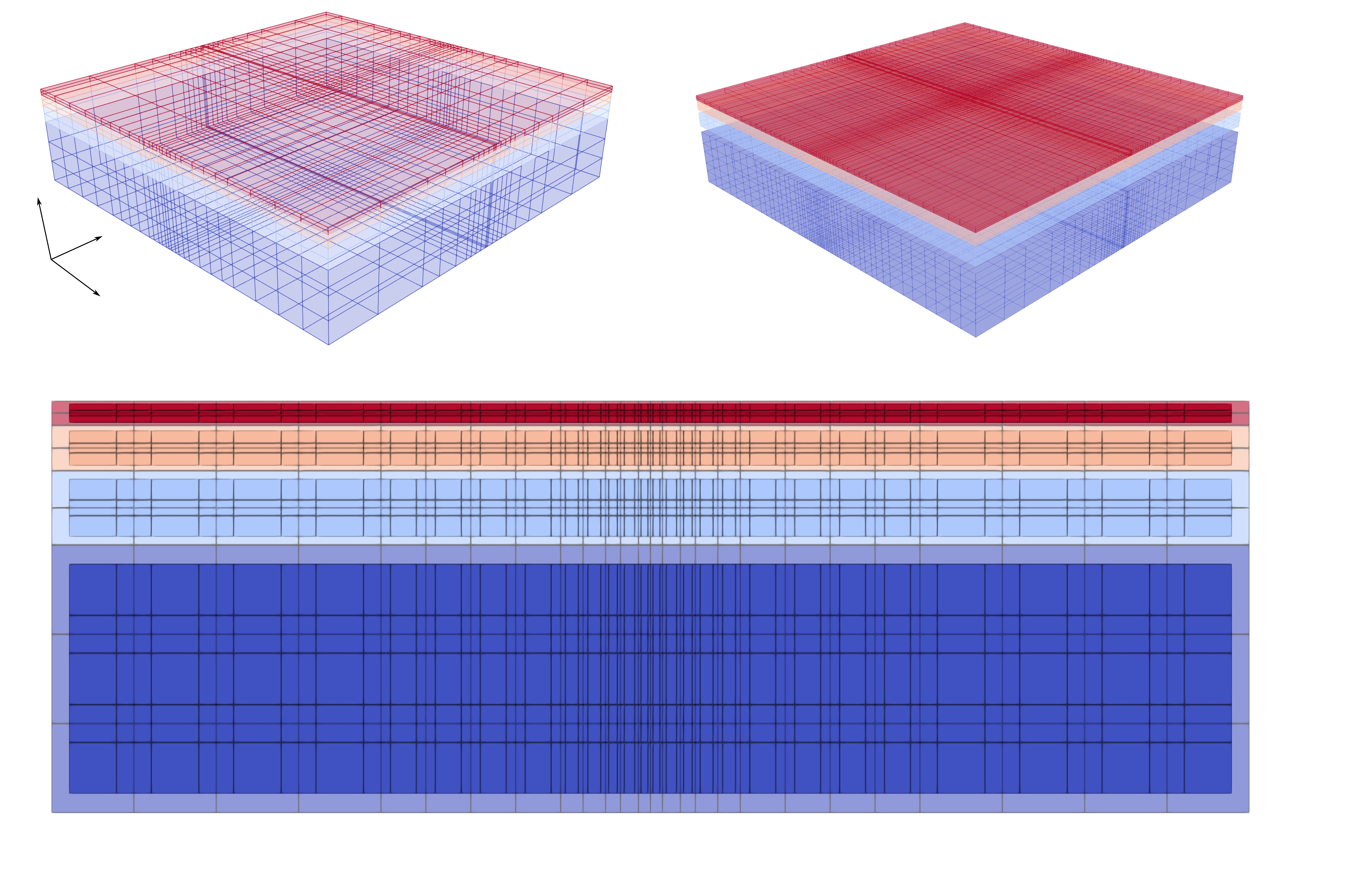}}%
		\put(0.01301556,0.64322938){\makebox(0,0)[lt]{\lineheight{1.25}\smash{\begin{tabular}[t]{l}(a)\end{tabular}}}}%
		\put(0.49432999,0.64321411){\makebox(0,0)[lt]{\lineheight{1.25}\smash{\begin{tabular}[t]{l}(b)\end{tabular}}}}%
		\put(0.01301556,0.37946462){\makebox(0,0)[lt]{\lineheight{1.25}\smash{\begin{tabular}[t]{l}(c)\end{tabular}}}}%
		\put(0.06988175,0.44100304){\makebox(0,0)[lt]{\lineheight{1.25}\smash{\begin{tabular}[t]{l}X\end{tabular}}}}%
		\put(0.07466369,0.47434018){\makebox(0,0)[lt]{\lineheight{1.25}\smash{\begin{tabular}[t]{l}Y\end{tabular}}}}%
		\put(0.0164374,0.51442818){\makebox(0,0)[lt]{\lineheight{1.25}\smash{\begin{tabular}[t]{l}Z\end{tabular}}}}%
		\put(0,0){\includegraphics[width=\unitlength,page=2]{02_gp_submesh_svg-tex.pdf}}%
		\put(0.530986,0.44100285){\makebox(0,0)[lt]{\lineheight{1.25}\smash{\begin{tabular}[t]{l}X\end{tabular}}}}%
		\put(0.53576795,0.47433949){\makebox(0,0)[lt]{\lineheight{1.25}\smash{\begin{tabular}[t]{l}Y\end{tabular}}}}%
		\put(0.47754165,0.51442754){\makebox(0,0)[lt]{\lineheight{1.25}\smash{\begin{tabular}[t]{l}Z\end{tabular}}}}%
		\put(0.0868731,0.01662412){\makebox(0,0)[lt]{\lineheight{1.25}\smash{\begin{tabular}[t]{l}X\end{tabular}}}}%
		\put(0.00264232,0.12203279){\makebox(0,0)[lt]{\lineheight{1.25}\smash{\begin{tabular}[t]{l}Z\end{tabular}}}}%
		\put(0.01606044,0.03224403){\makebox(0,0)[lt]{\lineheight{1.25}\smash{\begin{tabular}[t]{l}Y\end{tabular}}}}%
		\put(0,0){\includegraphics[width=\unitlength,page=3]{02_gp_submesh_svg-tex.pdf}}%
		\put(0.43942969,0.04124374){\makebox(0,0)[lt]{\lineheight{1.25}\smash{\begin{tabular}[t]{l}original finite element mesh (light colors)\end{tabular}}}}%
		\put(0,0){\includegraphics[width=\unitlength,page=4]{02_gp_submesh_svg-tex.pdf}}%
		\put(0.59212574,0.00659521){\makebox(0,0)[lt]{\lineheight{1.25}\smash{\begin{tabular}[t]{l}Gauss point sub-meshes for each layer (dark colors)\end{tabular}}}}%
	\end{picture}%
	\endgroup%
	\caption{(a) Perspective view of the original mesh, (b) perspective view of the Gauss point sub-meshes, and (c) side view with the Gauss point sub-meshes (for each layer) overlaid onto the original mesh, corresponding to the study in Section \ref{ssec:neovis}. The colours indicate different layers of the pavement structure.}
	\label{fig:submesh02}
\end{figure}

This algorithm can also be used for the update of internal variables $\boldsymbol{\alpha}$ in case the structure under investigation includes inelastic materials (as in the study reported in Section \ref{ssec:neovis}), if a sub-mesh is generated using the Gauss integration points as nodes \cite{anantheswar2025treatment}. The Gauss points of the main mesh are then represented as the nodes of the sub-mesh. Typically, the evolution of these internal variables is provided in terms of a rate equation involving the material time derivative, which is of the form
\begin{equation}\label{eq:alpha_mat_time_derv}
	\dot{\boldsymbol{\alpha}} =
	\left.\frac{\partial \boldsymbol{\alpha}}
	{\partial t}\right|_{\boldsymbol{X}} =
	\underbrace{\left.\frac{\partial \boldsymbol{\alpha}}
		{\partial t}\right|_{\boldsymbol{\chi}}}_{\text{Lagrangian phase}} +
	\underbrace{\text{Grad}\left( \boldsymbol{\alpha} \right) \cdot
		\boldsymbol{w}}_{\text{Eulerian phase}} \;\;.
\end{equation}
Using an operator split technique, the material time derivative in the Equation (\ref{eq:alpha_mat_time_derv}) can be structured into a Lagrangian phase and an Eulerian phase \cite{anantheswar2025treatment}. Further, the direct Gauss point advection (DGPA) algorithm utilized in this work for the history variable advection procedure \cite{anantheswar2025treatment} is also extended to be used with multilayered structures \cite{am3p_paper_atul}. For the non-uniform mesh simulated in Section \ref{ssec:neovis}, the original mesh and the Gauss point sub-mesh are visualized in Figure \ref{fig:submesh02}.

The most important aspect in simulating multilayered structures with different material models for the layers, is to isolate the Gauss point sub-meshes, such that the history variables of one material are not interpolated into those of another material. Figure \ref{fig:submesh02}(c) clearly shows this separation of the Gauss point sub-meshes by material. In this study, for the sake of convenience, the same kind of finite elements as in the original mesh (linear bricks with eight nodes) are used in the Gauss point sub-meshes. However, other kinds of finite elements may also be adopted. Additionally, the sub-mesh generation needs to be performed only once during the simulation. The generated Gauss point sub-mesh consists of 37440 nodes connected by 27846 finite elements. Thus, with this sub-mesh at hand, the history variables can be advected component wise through the mesh using the algorithm in Table \ref{table:disp_upd}.

\section{Application of Model Order Reduction techniques}\label{sec:MOR}
One particularly interesting aspect of the assembled system of equations in Equation (\ref{eq:soe}) is that it is relatively straightforward to apply projection-based model order reduction to reduce the complexity of solving the equation system, as shown in \cite{am3p_paper_jannick}. To create the reduced order model (ROM), it is assumed that a projection matrix $\underline{\boldsymbol{\Phi}} \in \mathbb{R}^{n \times m}$ can be found such that the displacements $\boldsymbol{\varphi} \in \mathbb{R}^{n \times 1}$ can be approximated as
\begin{equation}
    \boldsymbol{\varphi} \approx \underline{\boldsymbol{\Phi}} \, \boldsymbol{\varphi}_{red},
\end{equation}
where $\boldsymbol{\varphi}_{red} \in \mathbb{R}^{m \times 1}$ is the reduced displacement vector. Inserting this relation into Equation (\ref{eq:KTdyn}) and applying the Galerkin projection leads to the reduced system of equations
\begin{equation}
    \underline{\boldsymbol{\Phi}}^T \underline{\boldsymbol{K}}_{T, dyn}\underline{\boldsymbol{\Phi}} \,\Delta \boldsymbol{\varphi}_{red} = \underline{\boldsymbol{\Phi}}^T \mathbf{G}.
\end{equation}
The resulting number of equations has thereby been reduced from $n$ to $m$. If $m \ll n$, the computation time of solving the equation system is greatly reduced.

In this contribution, the proper orthogonal decomposition technique is used to compute the projection matrix $\underline{\boldsymbol{\Phi}}$. This is done by collecting $\ell$ solution states $\boldsymbol{\varphi}_i$ using the full order model (FOM) in a so-called snapshot matrix $\underline{\boldsymbol{D}} = [\boldsymbol{\varphi}_1 \, \dots \, \boldsymbol{\varphi}_{\ell}] \in \mathbb{R}^{n \times \ell}$. For example, the snapshots can be obtained from a single full order simulation, where the solution vector $\boldsymbol{\varphi}_i$ from all $n_t$ discrete time steps is aggregated into the snapshot matrix. After the snapshot matrix has been constructed, singular value decomposition is applied to the snapshot matrix
\begin{equation}
\underline{\boldsymbol{D}} = \underline{\boldsymbol{U}} \, \underline{\boldsymbol{\Sigma}} \, \underline{\boldsymbol{V}}^T,
\end{equation}
where $\underline{\boldsymbol{U}} \in \mathbb{R}^{n \times n}$ and $\underline{\boldsymbol{V}}^T \in \mathbb{R}^{\ell \times \ell}$ are matrices containing the left and right singular vectors, respectively. Additionally, $\underline{\boldsymbol{\Sigma}} \in \mathbb{R}^{n \times \ell}$ contains singular values in decreasing magnitude $\sigma_1 \ge \sigma_i \ge \sigma_{\ell}$ on its diagonal, where each singular value $\sigma_i$ indicates the importance of the corresponding singular vector $\boldsymbol{U}_i$, also called POD mode. Finally, the projection matrix $\underline{\boldsymbol{\Phi}}$ is obtained by taking the first $m$ vectors of the matrix of left singular vectors $\underline{\boldsymbol{U}}$, such that it reads $\underline{\boldsymbol{\Phi}} = [\boldsymbol{U}_1 \, \dots \, \boldsymbol{U}_m]$. The number of POD modes $m$ that should be included in the projection matrix can be determined by evaluating the singular values and finding the number of modes $m$, such that an energy tolerance is met according to
\begin{equation}
e = \frac{\sum^m_{j=1} \sigma_j}{\sum^l_{i=1} \sigma_i} \overset{!}{\ge} tol.
\end{equation}
Further, with the POD technique as a baseline, more advanced hyper-reduction techniques such as the Discrete Empirical Interpolation Method (DEIM) \cite{chaturantabut2010nonlinear} or the Energy-Conserving Sampling and Weighting method (ECSW) \cite{farhat2014dimensional} are typically used to further enhance computation speed. However, direct utilization of these techniques with the dynamic ALE formulation used in this contribution is not possible. Hyper-reduction techniques like DEIM and ECSW further reduce the problem by considering only the most significant finite elements of the structure and neglecting the rest. Doing so would lead to jumps or discontinuities in the solution fields, and more importantly, the fields of internal (history) variables, because the material constitutive law is not evaluated in the neglected finite elements. Then, this discontinuous representation of the field of internal (history) variables cannot be updated through the procedure outlined in Section \ref{sec:ALE_updates}. Therefore, in this manuscript, only POD based model order reduction is employed in the simulations.

\section{Numerical Studies}\label{sec:numex}
In this section, two numerical studies, which highlight the benefits of the proposed MORALE formulation are undertaken. Both studies involve transient simulations, and the Newmark time integration scheme is used with a time step size of 0.05 s and parameters $\beta = 0.25$, $\gamma = 0.5$ \cite{newmark1959method}. All simulations reported in this contribution adopt linear brick type finite elements with eight nodes and eight Gaussian integration points. For full order simulations, a direct sparse solver was used, whereas the dense reduced equation system was solved using a direct dense solver.

\subsection{Study 1: Basic example}\label{ssec:stvk}
\begin{figure}[h!]
	\centering
	\def\svgwidth{0.9\textwidth}
	\begingroup%
	\makeatletter%
	\providecommand\color[2][]{%
		\errmessage{(Inkscape) Color is used for the text in Inkscape, but the package 'color.sty' is not loaded}%
		\renewcommand\color[2][]{}%
	}%
	\providecommand\transparent[1]{%
		\errmessage{(Inkscape) Transparency is used (non-zero) for the text in Inkscape, but the package 'transparent.sty' is not loaded}%
		\renewcommand\transparent[1]{}%
	}%
	\providecommand\rotatebox[2]{#2}%
	\newcommand*\fsize{\dimexpr\f@size pt\relax}%
	\newcommand*\lineheight[1]{\fontsize{\fsize}{#1\fsize}\selectfont}%
	\ifx\svgwidth\undefined%
	\setlength{\unitlength}{1394.94999923bp}%
	\ifx\svgscale\undefined%
	\relax%
	\else%
	\setlength{\unitlength}{\unitlength * \real{\svgscale}}%
	\fi%
	\else%
	\setlength{\unitlength}{\svgwidth}%
	\fi%
	\global\let\svgwidth\undefined%
	\global\let\svgscale\undefined%
	\makeatother%
	\begin{picture}(1,0.78430869)%
		\lineheight{1}%
		\setlength\tabcolsep{0pt}%
		\put(0,0){\includegraphics[width=\unitlength,page=1]{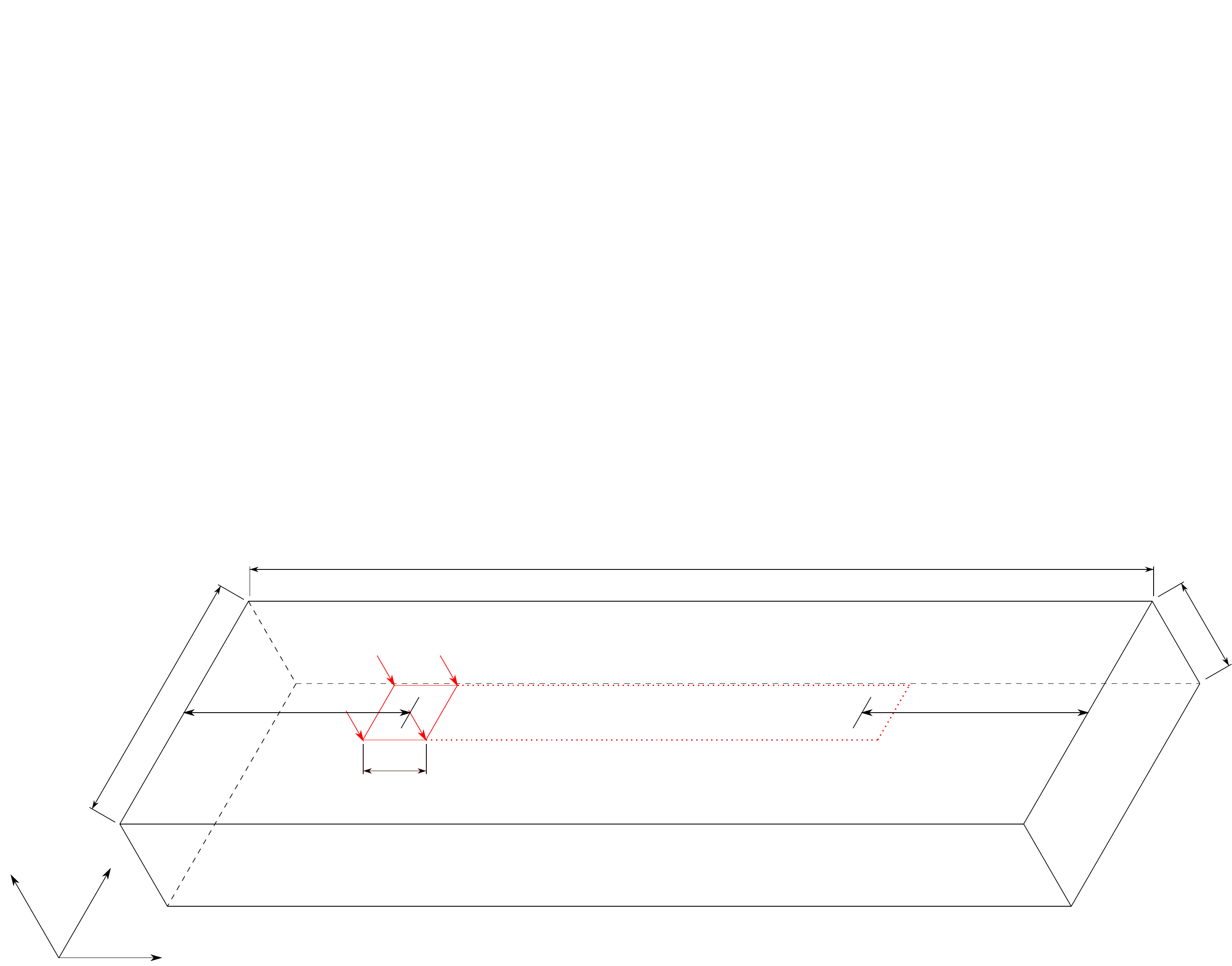}}%
		\put(0.09003008,0.36297105){\makebox(0,0)[lt]{\lineheight{1.25}\smash{\begin{tabular}[t]{l}(b)\end{tabular}}}}%
		\put(0,0){\includegraphics[width=\unitlength,page=2]{00_specimens_svg-tex.pdf}}%
		\put(0.07913251,0.2283252){\makebox(0,0)[lt]{\lineheight{1.25}\smash{\begin{tabular}[t]{l}160\end{tabular}}}}%
		\put(0.13945467,-0.00000002){\makebox(0,0)[lt]{\lineheight{1.25}\smash{\begin{tabular}[t]{l}X\end{tabular}}}}%
		\put(0.06760673,0.08001084){\makebox(0,0)[lt]{\lineheight{1.25}\smash{\begin{tabular}[t]{l}Y\end{tabular}}}}%
		\put(-0.00025763,0.07761794){\makebox(0,0)[lt]{\lineheight{1.25}\smash{\begin{tabular}[t]{l}Z\end{tabular}}}}%
		\put(0.32538329,0.2170872){\makebox(0,0)[lt]{\lineheight{1.25}\smash{\begin{tabular}[t]{l}p\end{tabular}}}}%
		\put(0.51994979,0.33954746){\makebox(0,0)[lt]{\lineheight{1.25}\smash{\begin{tabular}[t]{l}450\end{tabular}}}}%
		\put(0.24713531,0.18025075){\makebox(0,0)[lt]{\lineheight{1.25}\smash{\begin{tabular}[t]{l}80\end{tabular}}}}%
		\put(0.30759794,0.12397863){\makebox(0,0)[lt]{\lineheight{1.25}\smash{\begin{tabular}[t]{l}20\end{tabular}}}}%
		\put(0.78200678,0.18025075){\makebox(0,0)[lt]{\lineheight{1.25}\smash{\begin{tabular}[t]{l}80\end{tabular}}}}%
		\put(0.99219055,0.28710362){\makebox(0,0)[lt]{\lineheight{1.25}\smash{\begin{tabular}[t]{l}80\end{tabular}}}}%
		\put(0.42610009,0.15367134){\makebox(0,0)[lt]{\lineheight{1.25}\smash{\begin{tabular}[t]{l}load movement\end{tabular}}}}%
		\put(0,0){\includegraphics[width=\unitlength,page=3]{00_specimens_svg-tex.pdf}}%
		\put(0.09368047,0.76472634){\makebox(0,0)[lt]{\lineheight{1.25}\smash{\begin{tabular}[t]{l}(a)\end{tabular}}}}%
		\put(0.08278291,0.63657372){\makebox(0,0)[lt]{\lineheight{1.25}\smash{\begin{tabular}[t]{l}160\end{tabular}}}}%
		\put(0.29905311,0.75508309){\makebox(0,0)[lt]{\lineheight{1.25}\smash{\begin{tabular}[t]{l}160\end{tabular}}}}%
		\put(0.3017035,0.59898339){\makebox(0,0)[lt]{\lineheight{1.25}\smash{\begin{tabular}[t]{l}80\end{tabular}}}}%
		\put(0.23498985,0.54924725){\makebox(0,0)[lt]{\lineheight{1.25}\smash{\begin{tabular}[t]{l}20\end{tabular}}}}%
		\put(0.1813069,0.64725178){\makebox(0,0)[lt]{\lineheight{1.25}\smash{\begin{tabular}[t]{l}20\end{tabular}}}}%
		\put(0.46425719,0.70525206){\makebox(0,0)[lt]{\lineheight{1.25}\smash{\begin{tabular}[t]{l}80\end{tabular}}}}%
		\put(0.54942191,0.58554693){\makebox(0,0)[lt]{\lineheight{1.25}\smash{\begin{tabular}[t]{l}guiding velocity $\boldsymbol{w}$\end{tabular}}}}%
		\put(0.14305997,0.41810759){\makebox(0,0)[lt]{\lineheight{1.25}\smash{\begin{tabular}[t]{l}X\end{tabular}}}}%
		\put(0.07125715,0.49166644){\makebox(0,0)[lt]{\lineheight{1.25}\smash{\begin{tabular}[t]{l}Y\end{tabular}}}}%
		\put(0.00334768,0.49572557){\makebox(0,0)[lt]{\lineheight{1.25}\smash{\begin{tabular}[t]{l}Z\end{tabular}}}}%
		\put(0.2470616,0.63714164){\makebox(0,0)[lt]{\lineheight{1.25}\smash{\begin{tabular}[t]{l}p\end{tabular}}}}%
		\put(0.79326221,0.58611858){\makebox(0,0)[lt]{\lineheight{1.25}\smash{\begin{tabular}[t]{l}all dimensions in [m]\end{tabular}}}}%
		\put(0.79331883,0.53916357){\makebox(0,0)[lt]{\lineheight{1.25}\smash{\begin{tabular}[t]{l}all surfaces except top surface\\are fully restrained\end{tabular}}}}%
	\end{picture}%
	\endgroup%
	\caption{Specimens used in (a) ALE / MORALE simulations and (b) Lagrangian simulation of Section \ref{ssec:stvk}.}
	\label{fig:specimens00}
\end{figure}

In this study, a comparison of results from conventional Lagrangian, ALE and MORALE formulations is made by performing simulations on a simple uniform mesh. For this study, a hyperelastic St. Venant Kirchhoff material model is used. The specimen geometry is depicted in Figure \ref{fig:specimens00}.

\begin{figure}[h!]
	\centering
	\def\svgwidth{\textwidth}
	\begingroup%
	\makeatletter%
	\providecommand\color[2][]{%
		\errmessage{(Inkscape) Color is used for the text in Inkscape, but the package 'color.sty' is not loaded}%
		\renewcommand\color[2][]{}%
	}%
	\providecommand\transparent[1]{%
		\errmessage{(Inkscape) Transparency is used (non-zero) for the text in Inkscape, but the package 'transparent.sty' is not loaded}%
		\renewcommand\transparent[1]{}%
	}%
	\providecommand\rotatebox[2]{#2}%
	\newcommand*\fsize{\dimexpr\f@size pt\relax}%
	\newcommand*\lineheight[1]{\fontsize{\fsize}{#1\fsize}\selectfont}%
	\ifx\svgwidth\undefined%
	\setlength{\unitlength}{590.38244629bp}%
	\ifx\svgscale\undefined%
	\relax%
	\else%
	\setlength{\unitlength}{\unitlength * \real{\svgscale}}%
	\fi%
	\else%
	\setlength{\unitlength}{\svgwidth}%
	\fi%
	\global\let\svgwidth\undefined%
	\global\let\svgscale\undefined%
	\makeatother%
	\begin{picture}(1,0.52033211)%
		\lineheight{1}%
		\setlength\tabcolsep{0pt}%
		\put(0,0){\includegraphics[width=\unitlength,page=1]{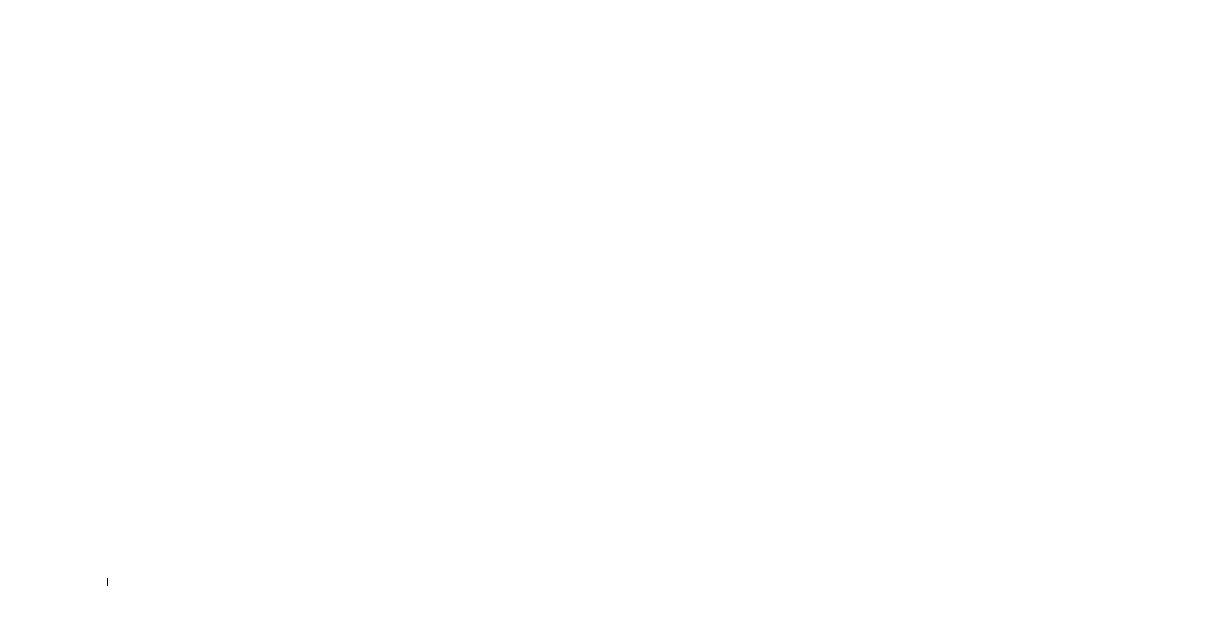}}%
		\put(0.08740229,0.02540134){\makebox(0,0)[t]{\lineheight{0}\smash{\begin{tabular}[t]{c}0\end{tabular}}}}%
		\put(0,0){\includegraphics[width=\unitlength,page=2]{ALEPOD_newnew_modified_svg-tex.pdf}}%
		\put(0.18927474,0.02540134){\makebox(0,0)[t]{\lineheight{0}\smash{\begin{tabular}[t]{c}2\end{tabular}}}}%
		\put(0,0){\includegraphics[width=\unitlength,page=3]{ALEPOD_newnew_modified_svg-tex.pdf}}%
		\put(0.29114722,0.02540134){\makebox(0,0)[t]{\lineheight{0}\smash{\begin{tabular}[t]{c}4\end{tabular}}}}%
		\put(0,0){\includegraphics[width=\unitlength,page=4]{ALEPOD_newnew_modified_svg-tex.pdf}}%
		\put(0.3930197,0.02540134){\makebox(0,0)[t]{\lineheight{0}\smash{\begin{tabular}[t]{c}6\end{tabular}}}}%
		\put(0,0){\includegraphics[width=\unitlength,page=5]{ALEPOD_newnew_modified_svg-tex.pdf}}%
		\put(0.49489215,0.02540134){\makebox(0,0)[t]{\lineheight{0}\smash{\begin{tabular}[t]{c}8\end{tabular}}}}%
		\put(0,0){\includegraphics[width=\unitlength,page=6]{ALEPOD_newnew_modified_svg-tex.pdf}}%
		\put(0.5967646,0.02540134){\makebox(0,0)[t]{\lineheight{0}\smash{\begin{tabular}[t]{c}10\end{tabular}}}}%
		\put(0,0){\includegraphics[width=\unitlength,page=7]{ALEPOD_newnew_modified_svg-tex.pdf}}%
		\put(0.69863705,0.02540134){\makebox(0,0)[t]{\lineheight{0}\smash{\begin{tabular}[t]{c}12\end{tabular}}}}%
		\put(0,0){\includegraphics[width=\unitlength,page=8]{ALEPOD_newnew_modified_svg-tex.pdf}}%
		\put(0.8005095,0.02540134){\makebox(0,0)[t]{\lineheight{0}\smash{\begin{tabular}[t]{c}14\end{tabular}}}}%
		\put(0,0){\includegraphics[width=\unitlength,page=9]{ALEPOD_newnew_modified_svg-tex.pdf}}%
		\put(0.90238201,0.02540134){\makebox(0,0)[t]{\lineheight{0}\smash{\begin{tabular}[t]{c}16\end{tabular}}}}%
		\put(0.51271982,0.00223308){\makebox(0,0)[t]{\lineheight{0}\smash{\begin{tabular}[t]{c}Time [s]\end{tabular}}}}%
		\put(0,0){\includegraphics[width=\unitlength,page=10]{ALEPOD_newnew_modified_svg-tex.pdf}}%
		\put(0.07554557,0.07514636){\makebox(0,0)[rt]{\lineheight{0}\smash{\begin{tabular}[t]{r}-0.35\end{tabular}}}}%
		\put(0,0){\includegraphics[width=\unitlength,page=11]{ALEPOD_newnew_modified_svg-tex.pdf}}%
		\put(0.07554557,0.13467933){\makebox(0,0)[rt]{\lineheight{0}\smash{\begin{tabular}[t]{r}-0.30\end{tabular}}}}%
		\put(0,0){\includegraphics[width=\unitlength,page=12]{ALEPOD_newnew_modified_svg-tex.pdf}}%
		\put(0.07554557,0.19421234){\makebox(0,0)[rt]{\lineheight{0}\smash{\begin{tabular}[t]{r}-0.25\end{tabular}}}}%
		\put(0,0){\includegraphics[width=\unitlength,page=13]{ALEPOD_newnew_modified_svg-tex.pdf}}%
		\put(0.07554557,0.25374533){\makebox(0,0)[rt]{\lineheight{0}\smash{\begin{tabular}[t]{r}-0.20\end{tabular}}}}%
		\put(0,0){\includegraphics[width=\unitlength,page=14]{ALEPOD_newnew_modified_svg-tex.pdf}}%
		\put(0.07554557,0.31327832){\makebox(0,0)[rt]{\lineheight{0}\smash{\begin{tabular}[t]{r}-0.15\end{tabular}}}}%
		\put(0,0){\includegraphics[width=\unitlength,page=15]{ALEPOD_newnew_modified_svg-tex.pdf}}%
		\put(0.07554557,0.37281131){\makebox(0,0)[rt]{\lineheight{0}\smash{\begin{tabular}[t]{r}-0.10\end{tabular}}}}%
		\put(0,0){\includegraphics[width=\unitlength,page=16]{ALEPOD_newnew_modified_svg-tex.pdf}}%
		\put(0.07554557,0.43234431){\makebox(0,0)[rt]{\lineheight{0}\smash{\begin{tabular}[t]{r}-0.05\end{tabular}}}}%
		\put(0,0){\includegraphics[width=\unitlength,page=17]{ALEPOD_newnew_modified_svg-tex.pdf}}%
		\put(0.07554557,0.4918773){\makebox(0,0)[rt]{\lineheight{0}\smash{\begin{tabular}[t]{r}0.00\end{tabular}}}}%
		\put(0.01286904,0.28489151){\rotatebox{90}{\makebox(0,0)[t]{\lineheight{0}\smash{\begin{tabular}[t]{c}Displacement $u_{z}$ [m] under load\end{tabular}}}}}%
		\put(0,0){\includegraphics[width=\unitlength,page=18]{ALEPOD_newnew_modified_svg-tex.pdf}}%
		\put(0.33987714,0.14348156){\makebox(0,0)[lt]{\lineheight{0}\smash{\begin{tabular}[t]{l}Lagrangian\end{tabular}}}}%
		\put(0,0){\includegraphics[width=\unitlength,page=19]{ALEPOD_newnew_modified_svg-tex.pdf}}%
		\put(0.33987714,0.11861954){\makebox(0,0)[lt]{\lineheight{0}\smash{\begin{tabular}[t]{l}POD-Lagrangian\end{tabular}}}}%
		\put(0,0){\includegraphics[width=\unitlength,page=20]{ALEPOD_newnew_modified_svg-tex.pdf}}%
		\put(0.33987714,0.09375747){\makebox(0,0)[lt]{\lineheight{0}\smash{\begin{tabular}[t]{l}ALE\end{tabular}}}}%
		\put(0,0){\includegraphics[width=\unitlength,page=21]{ALEPOD_newnew_modified_svg-tex.pdf}}%
		\put(0.33987714,0.0688954){\makebox(0,0)[lt]{\lineheight{0}\smash{\begin{tabular}[t]{l}MORALE\end{tabular}}}}%
		\put(0,0){\includegraphics[width=\unitlength,page=22]{ALEPOD_newnew_modified_svg-tex.pdf}}%
		\put(0.94989408,0.06503536){\makebox(0,0)[lt]{\lineheight{0}\smash{\begin{tabular}[t]{l}0.0\end{tabular}}}}%
		\put(0,0){\includegraphics[width=\unitlength,page=23]{ALEPOD_newnew_modified_svg-tex.pdf}}%
		\put(0.94989408,0.15040374){\makebox(0,0)[lt]{\lineheight{0}\smash{\begin{tabular}[t]{l}0.2\end{tabular}}}}%
		\put(0,0){\includegraphics[width=\unitlength,page=24]{ALEPOD_newnew_modified_svg-tex.pdf}}%
		\put(0.94989408,0.23577212){\makebox(0,0)[lt]{\lineheight{0}\smash{\begin{tabular}[t]{l}0.4\end{tabular}}}}%
		\put(0,0){\includegraphics[width=\unitlength,page=25]{ALEPOD_newnew_modified_svg-tex.pdf}}%
		\put(0.94989408,0.32114053){\makebox(0,0)[lt]{\lineheight{0}\smash{\begin{tabular}[t]{l}0.6\end{tabular}}}}%
		\put(0,0){\includegraphics[width=\unitlength,page=26]{ALEPOD_newnew_modified_svg-tex.pdf}}%
		\put(0.94989408,0.40650891){\makebox(0,0)[lt]{\lineheight{0}\smash{\begin{tabular}[t]{l}0.8\end{tabular}}}}%
		\put(0,0){\includegraphics[width=\unitlength,page=27]{ALEPOD_newnew_modified_svg-tex.pdf}}%
		\put(0.94989408,0.4918773){\makebox(0,0)[lt]{\lineheight{0}\smash{\begin{tabular}[t]{l}1.0\end{tabular}}}}%
		\put(0.99647668,0.28489147){\rotatebox{90}{\makebox(0,0)[t]{\lineheight{0}\smash{\begin{tabular}[t]{c}Amplitude scaling factor [-]\end{tabular}}}}}%
		\put(0,0){\includegraphics[width=\unitlength,page=28]{ALEPOD_newnew_modified_svg-tex.pdf}}%
		\put(0.53544307,0.09375747){\makebox(0,0)[lt]{\lineheight{0}\smash{\begin{tabular}[t]{l}Pressure load\end{tabular}}}}%
		\put(0,0){\includegraphics[width=\unitlength,page=29]{ALEPOD_newnew_modified_svg-tex.pdf}}%
		\put(0.53544307,0.0688954){\makebox(0,0)[lt]{\lineheight{0}\smash{\begin{tabular}[t]{l}Load movement / Guiding velocity\end{tabular}}}}%
	\end{picture}%
	\endgroup%
    \caption{Displacement $u_z$ under the center of the load for the study in Section \ref{ssec:stvk}, plotted over time.}
	\label{fig:dispvstime00}
\end{figure}

\begin{figure}[h!]
	\centering
	\def\svgwidth{\textwidth}
	\begingroup%
	\makeatletter%
	\providecommand\color[2][]{%
		\errmessage{(Inkscape) Color is used for the text in Inkscape, but the package 'color.sty' is not loaded}%
		\renewcommand\color[2][]{}%
	}%
	\providecommand\transparent[1]{%
		\errmessage{(Inkscape) Transparency is used (non-zero) for the text in Inkscape, but the package 'transparent.sty' is not loaded}%
		\renewcommand\transparent[1]{}%
	}%
	\providecommand\rotatebox[2]{#2}%
	\newcommand*\fsize{\dimexpr\f@size pt\relax}%
	\newcommand*\lineheight[1]{\fontsize{\fsize}{#1\fsize}\selectfont}%
	\ifx\svgwidth\undefined%
	\setlength{\unitlength}{554.20080566bp}%
	\ifx\svgscale\undefined%
	\relax%
	\else%
	\setlength{\unitlength}{\unitlength * \real{\svgscale}}%
	\fi%
	\else%
	\setlength{\unitlength}{\svgwidth}%
	\fi%
	\global\let\svgwidth\undefined%
	\global\let\svgscale\undefined%
	\makeatother%
	\begin{picture}(1,0.5543026)%
		\lineheight{1}%
		\setlength\tabcolsep{0pt}%
		\put(0,0){\includegraphics[width=\unitlength,page=1]{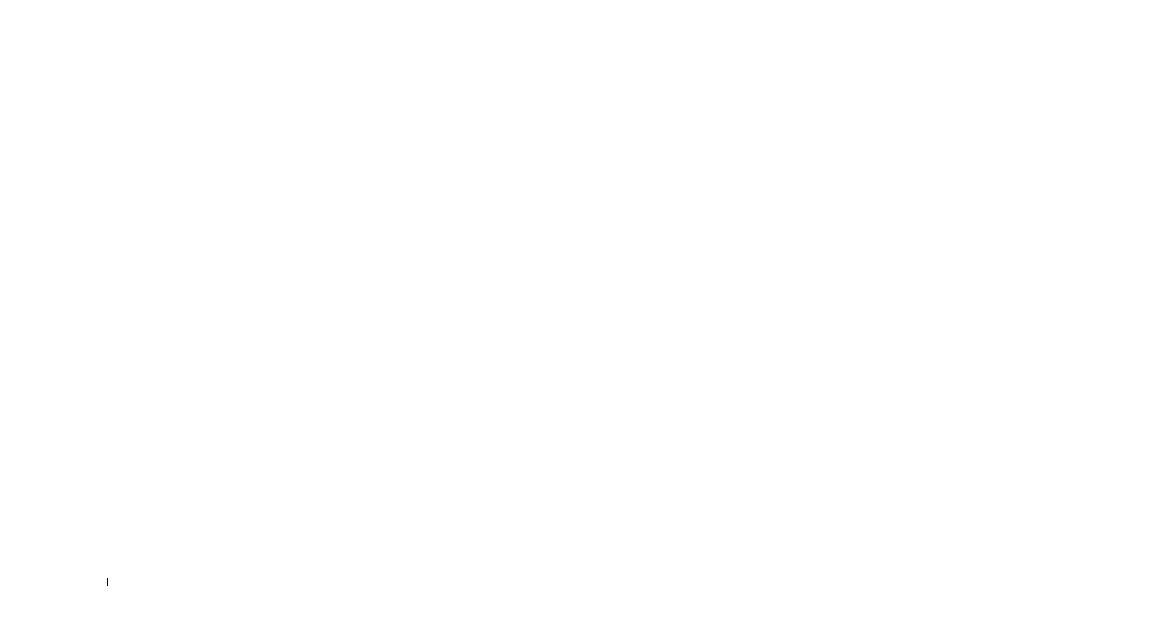}}%
		\put(0.09310845,0.02705969){\makebox(0,0)[t]{\lineheight{1.25}\smash{\begin{tabular}[t]{c}0\end{tabular}}}}%
		\put(0,0){\includegraphics[width=\unitlength,page=2]{ALEPOD_modified_svg-tex.pdf}}%
		\put(0.20163176,0.02705969){\makebox(0,0)[t]{\lineheight{1.25}\smash{\begin{tabular}[t]{c}2\end{tabular}}}}%
		\put(0,0){\includegraphics[width=\unitlength,page=3]{ALEPOD_modified_svg-tex.pdf}}%
		\put(0.3101551,0.02705969){\makebox(0,0)[t]{\lineheight{1.25}\smash{\begin{tabular}[t]{c}4\end{tabular}}}}%
		\put(0,0){\includegraphics[width=\unitlength,page=4]{ALEPOD_modified_svg-tex.pdf}}%
		\put(0.41867844,0.02705969){\makebox(0,0)[t]{\lineheight{1.25}\smash{\begin{tabular}[t]{c}6\end{tabular}}}}%
		\put(0,0){\includegraphics[width=\unitlength,page=5]{ALEPOD_modified_svg-tex.pdf}}%
		\put(0.52720176,0.02705969){\makebox(0,0)[t]{\lineheight{1.25}\smash{\begin{tabular}[t]{c}8\end{tabular}}}}%
		\put(0,0){\includegraphics[width=\unitlength,page=6]{ALEPOD_modified_svg-tex.pdf}}%
		\put(0.63572507,0.02705969){\makebox(0,0)[t]{\lineheight{1.25}\smash{\begin{tabular}[t]{c}10\end{tabular}}}}%
		\put(0,0){\includegraphics[width=\unitlength,page=7]{ALEPOD_modified_svg-tex.pdf}}%
		\put(0.74424838,0.02705969){\makebox(0,0)[t]{\lineheight{1.25}\smash{\begin{tabular}[t]{c}12\end{tabular}}}}%
		\put(0,0){\includegraphics[width=\unitlength,page=8]{ALEPOD_modified_svg-tex.pdf}}%
		\put(0.85277169,0.02705969){\makebox(0,0)[t]{\lineheight{1.25}\smash{\begin{tabular}[t]{c}14\end{tabular}}}}%
		\put(0,0){\includegraphics[width=\unitlength,page=9]{ALEPOD_modified_svg-tex.pdf}}%
		\put(0.96129506,0.02705969){\makebox(0,0)[t]{\lineheight{1.25}\smash{\begin{tabular}[t]{c}16\end{tabular}}}}%
		\put(0.54619333,0.00237887){\makebox(0,0)[t]{\lineheight{1.25}\smash{\begin{tabular}[t]{c}Time [s]\end{tabular}}}}%
		\put(0,0){\includegraphics[width=\unitlength,page=10]{ALEPOD_modified_svg-tex.pdf}}%
		\put(0.08047765,0.07495069){\makebox(0,0)[rt]{\lineheight{1.25}\smash{\begin{tabular}[t]{r}-0.35\end{tabular}}}}%
		\put(0,0){\includegraphics[width=\unitlength,page=11]{ALEPOD_modified_svg-tex.pdf}}%
		\put(0.08047765,0.13909918){\makebox(0,0)[rt]{\lineheight{1.25}\smash{\begin{tabular}[t]{r}-0.30\end{tabular}}}}%
		\put(0,0){\includegraphics[width=\unitlength,page=12]{ALEPOD_modified_svg-tex.pdf}}%
		\put(0.08047765,0.20324765){\makebox(0,0)[rt]{\lineheight{1.25}\smash{\begin{tabular}[t]{r}-0.25\end{tabular}}}}%
		\put(0,0){\includegraphics[width=\unitlength,page=13]{ALEPOD_modified_svg-tex.pdf}}%
		\put(0.08047765,0.26739615){\makebox(0,0)[rt]{\lineheight{1.25}\smash{\begin{tabular}[t]{r}-0.20\end{tabular}}}}%
		\put(0,0){\includegraphics[width=\unitlength,page=14]{ALEPOD_modified_svg-tex.pdf}}%
		\put(0.08047765,0.33154464){\makebox(0,0)[rt]{\lineheight{1.25}\smash{\begin{tabular}[t]{r}-0.15\end{tabular}}}}%
		\put(0,0){\includegraphics[width=\unitlength,page=15]{ALEPOD_modified_svg-tex.pdf}}%
		\put(0.08047765,0.39569311){\makebox(0,0)[rt]{\lineheight{1.25}\smash{\begin{tabular}[t]{r}-0.10\end{tabular}}}}%
		\put(0,0){\includegraphics[width=\unitlength,page=16]{ALEPOD_modified_svg-tex.pdf}}%
		\put(0.08047765,0.4598416){\makebox(0,0)[rt]{\lineheight{1.25}\smash{\begin{tabular}[t]{r}-0.05\end{tabular}}}}%
		\put(0,0){\includegraphics[width=\unitlength,page=17]{ALEPOD_modified_svg-tex.pdf}}%
		\put(0.08047765,0.52399008){\makebox(0,0)[rt]{\lineheight{1.25}\smash{\begin{tabular}[t]{r}0.00\end{tabular}}}}%
		\put(0.01370921,0.30349098){\rotatebox{90}{\makebox(0,0)[t]{\lineheight{1.25}\smash{\begin{tabular}[t]{c}Displacement $u_{z}$ [m] under load\end{tabular}}}}}%
		\put(0,0){\includegraphics[width=\unitlength,page=18]{ALEPOD_modified_svg-tex.pdf}}%
		\put(0.48896481,0.31677869){\makebox(0,0)[lt]{\lineheight{1.25}\smash{\begin{tabular}[t]{l}FOM $w_{max} = 25$ m/s, $p_{max} = 50$ Pa\end{tabular}}}}%
		\put(0,0){\includegraphics[width=\unitlength,page=19]{ALEPOD_modified_svg-tex.pdf}}%
		\put(0.48896481,0.28979165){\makebox(0,0)[lt]{\lineheight{1.25}\smash{\begin{tabular}[t]{l}POD $w_{max} = 25$ m/s, $p_{max} = 50$ Pa\end{tabular}}}}%
		\put(0,0){\includegraphics[width=\unitlength,page=20]{ALEPOD_modified_svg-tex.pdf}}%
		\put(0.48896481,0.26280458){\makebox(0,0)[lt]{\lineheight{1.25}\smash{\begin{tabular}[t]{l}FOM $w_{max} = 20$ m/s, $p_{max} = 50$ Pa\end{tabular}}}}%
		\put(0,0){\includegraphics[width=\unitlength,page=21]{ALEPOD_modified_svg-tex.pdf}}%
		\put(0.48896481,0.23581752){\makebox(0,0)[lt]{\lineheight{1.25}\smash{\begin{tabular}[t]{l}POD $w_{max} = 20$ m/s, $p_{max} = 50$ Pa\end{tabular}}}}%
		\put(0,0){\includegraphics[width=\unitlength,page=22]{ALEPOD_modified_svg-tex.pdf}}%
		\put(0.48896481,0.20883045){\makebox(0,0)[lt]{\lineheight{1.25}\smash{\begin{tabular}[t]{l}FOM $w_{max} = 30$ m/s, $p_{max} = 50$ Pa\end{tabular}}}}%
		\put(0,0){\includegraphics[width=\unitlength,page=23]{ALEPOD_modified_svg-tex.pdf}}%
		\put(0.48896481,0.18184339){\makebox(0,0)[lt]{\lineheight{1.25}\smash{\begin{tabular}[t]{l}POD $w_{max} = 30$ m/s, $p_{max} = 50$ Pa\end{tabular}}}}%
		\put(0,0){\includegraphics[width=\unitlength,page=24]{ALEPOD_modified_svg-tex.pdf}}%
		\put(0.48896481,0.15485635){\makebox(0,0)[lt]{\lineheight{1.25}\smash{\begin{tabular}[t]{l}FOM $w_{max} = 20$ m/s, $p_{max} = 25$ Pa\end{tabular}}}}%
		\put(0,0){\includegraphics[width=\unitlength,page=25]{ALEPOD_modified_svg-tex.pdf}}%
		\put(0.48896481,0.12786926){\makebox(0,0)[lt]{\lineheight{1.25}\smash{\begin{tabular}[t]{l}POD $w_{max} = 20$ m/s, $p_{max} = 25$ Pa\end{tabular}}}}%
		\put(0,0){\includegraphics[width=\unitlength,page=26]{ALEPOD_modified_svg-tex.pdf}}%
		\put(0.48896481,0.10088222){\makebox(0,0)[lt]{\lineheight{1.25}\smash{\begin{tabular}[t]{l}FOM $w_{max} = 30$ m/s, $p_{max} = 25$ Pa\end{tabular}}}}%
		\put(0,0){\includegraphics[width=\unitlength,page=27]{ALEPOD_modified_svg-tex.pdf}}%
		\put(0.48896481,0.07389518){\makebox(0,0)[lt]{\lineheight{1.25}\smash{\begin{tabular}[t]{l}POD $w_{max} = 30$ m/s, $p_{max} = 25$ Pa\end{tabular}}}}%
	\end{picture}%
	\endgroup%
	\caption{Predictive qualities of the reduced order model created only from the snapshots of the full order ALE simulation with maximum guiding velocity $w_{max} = 25$ m/s and load $p = 50$ Pa, when subjected to varying guiding velocities $w_{max}$ and loads $p$. Here, FOM stands for `Full Order Model'.}
	\label{fig:POD_performance}
\end{figure}

\begin{figure}[h!]
	\centering
	\def\svgwidth{0.9\textwidth}
	\begingroup%
	\makeatletter%
	\providecommand\color[2][]{%
		\errmessage{(Inkscape) Color is used for the text in Inkscape, but the package 'color.sty' is not loaded}%
		\renewcommand\color[2][]{}%
	}%
	\providecommand\transparent[1]{%
		\errmessage{(Inkscape) Transparency is used (non-zero) for the text in Inkscape, but the package 'transparent.sty' is not loaded}%
		\renewcommand\transparent[1]{}%
	}%
	\providecommand\rotatebox[2]{#2}%
	\newcommand*\fsize{\dimexpr\f@size pt\relax}%
	\newcommand*\lineheight[1]{\fontsize{\fsize}{#1\fsize}\selectfont}%
	\ifx\svgwidth\undefined%
	\setlength{\unitlength}{460.79998779bp}%
	\ifx\svgscale\undefined%
	\relax%
	\else%
	\setlength{\unitlength}{\unitlength * \real{\svgscale}}%
	\fi%
	\else%
	\setlength{\unitlength}{\svgwidth}%
	\fi%
	\global\let\svgwidth\undefined%
	\global\let\svgscale\undefined%
	\makeatother%
	\begin{picture}(1,0.75000003)%
		\lineheight{1}%
		\setlength\tabcolsep{0pt}%
		\put(0,0){\includegraphics[width=\unitlength,page=1]{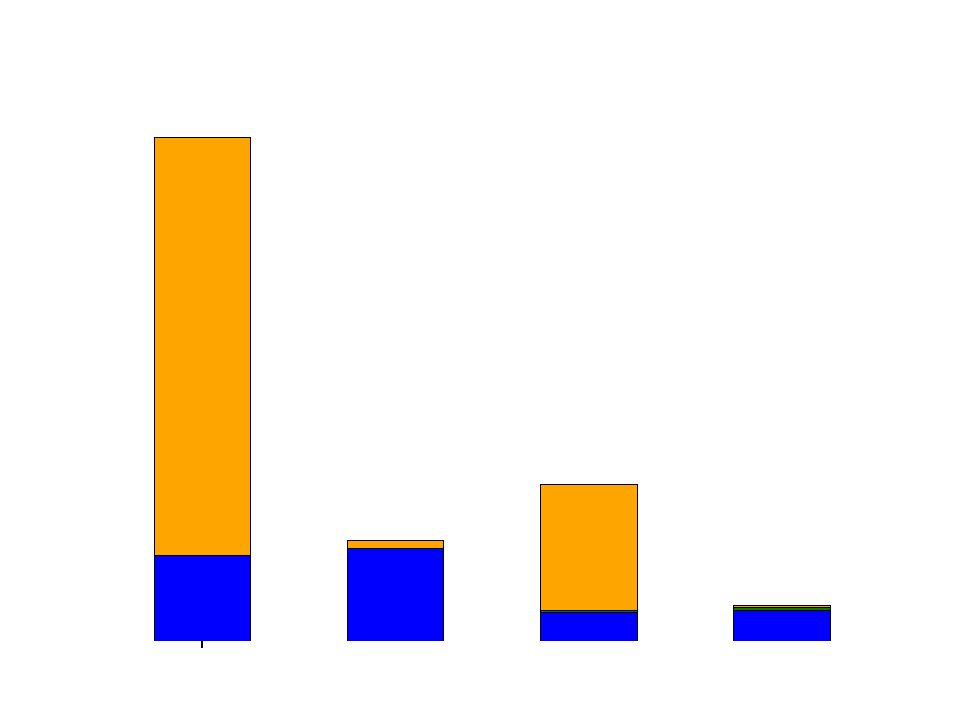}}%
		\put(0.21055196,0.05081939){\makebox(0,0)[t]{\lineheight{1.25}\smash{\begin{tabular}[t]{c}Lagrange\end{tabular}}}}%
		\put(0,0){\includegraphics[width=\unitlength,page=2]{1_timings_svg-tex.pdf}}%
		\put(0.41185066,0.05081939){\makebox(0,0)[t]{\lineheight{1.25}\smash{\begin{tabular}[t]{c}POD-Lagrange\end{tabular}}}}%
		\put(0,0){\includegraphics[width=\unitlength,page=3]{1_timings_svg-tex.pdf}}%
		\put(0.61314935,0.05081939){\makebox(0,0)[t]{\lineheight{1.25}\smash{\begin{tabular}[t]{c}ALE\end{tabular}}}}%
		\put(0,0){\includegraphics[width=\unitlength,page=4]{1_timings_svg-tex.pdf}}%
		\put(0.78539023,0.05081939){\makebox(0,0)[lt]{\lineheight{1.25}\smash{\begin{tabular}[t]{l}MORALE\end{tabular}}}}%
		\put(0.78829789,0.02740019){\makebox(0,0)[lt]{\lineheight{1.25}\smash{\begin{tabular}[t]{l}($m$ = 73)\end{tabular}}}}%
		\put(0,0){\includegraphics[width=\unitlength,page=5]{1_timings_svg-tex.pdf}}%
		\put(0.10980903,0.0742552){\makebox(0,0)[rt]{\lineheight{1.25}\smash{\begin{tabular}[t]{r}0.0\end{tabular}}}}%
		\put(0,0){\includegraphics[width=\unitlength,page=6]{1_timings_svg-tex.pdf}}%
		\put(0.10980903,0.17925521){\makebox(0,0)[rt]{\lineheight{1.25}\smash{\begin{tabular}[t]{r}0.2\end{tabular}}}}%
		\put(0,0){\includegraphics[width=\unitlength,page=7]{1_timings_svg-tex.pdf}}%
		\put(0.10980903,0.28425519){\makebox(0,0)[rt]{\lineheight{1.25}\smash{\begin{tabular}[t]{r}0.4\end{tabular}}}}%
		\put(0,0){\includegraphics[width=\unitlength,page=8]{1_timings_svg-tex.pdf}}%
		\put(0.10980903,0.3892552){\makebox(0,0)[rt]{\lineheight{1.25}\smash{\begin{tabular}[t]{r}0.6\end{tabular}}}}%
		\put(0,0){\includegraphics[width=\unitlength,page=9]{1_timings_svg-tex.pdf}}%
		\put(0.10980903,0.49425519){\makebox(0,0)[rt]{\lineheight{1.25}\smash{\begin{tabular}[t]{r}0.8\end{tabular}}}}%
		\put(0,0){\includegraphics[width=\unitlength,page=10]{1_timings_svg-tex.pdf}}%
		\put(0.10980903,0.5992552){\makebox(0,0)[rt]{\lineheight{1.25}\smash{\begin{tabular}[t]{r}1.0\end{tabular}}}}%
		\put(0.06210327,0.37125002){\rotatebox{90}{\makebox(0,0)[t]{\lineheight{1.25}\smash{\begin{tabular}[t]{c}Relative simulation time [-]\end{tabular}}}}}%
		\put(0,0){\includegraphics[width=\unitlength,page=11]{1_timings_svg-tex.pdf}}%
		\put(0.43640747,0.62397911){\makebox(0,0)[lt]{\lineheight{1.25}\smash{\begin{tabular}[t]{l}Computing (reduced) material quantities\end{tabular}}}}%
		\put(0,0){\includegraphics[width=\unitlength,page=12]{1_timings_svg-tex.pdf}}%
		\put(0.43640747,0.59212554){\makebox(0,0)[lt]{\lineheight{1.25}\smash{\begin{tabular}[t]{l}Computing (reduced) ALE quantities\end{tabular}}}}%
		\put(0,0){\includegraphics[width=\unitlength,page=13]{1_timings_svg-tex.pdf}}%
		\put(0.43640747,0.56027197){\makebox(0,0)[lt]{\lineheight{1.25}\smash{\begin{tabular}[t]{l}Solving (reduced) system of equations\end{tabular}}}}%
		\put(0,0){\includegraphics[width=\unitlength,page=14]{1_timings_svg-tex.pdf}}%
		\put(0.43640747,0.5284184){\makebox(0,0)[lt]{\lineheight{1.25}\smash{\begin{tabular}[t]{l}POD basis computation\end{tabular}}}}%
	\end{picture}%
	\endgroup%
	\caption{Comparison of times required for simulations in Section \ref{ssec:stvk}.}
	\label{fig:runtimes00}
\end{figure}

The boundary conditions used in Lagrangian, ALE and MORALE simulations are such that all surfaces of the simulated specimen, except the top surface (where pressure load is applied), are fully restrained. The material parameters used are $\lambda = 2500$ Pa, $\mu = 1250$ Pa (Lamé parameters) and density $\hat{\rho} = 10.0$ kg/m$^3$. In the Lagrangian simulation, the pressure load p = 50 Pa is first applied as a linear ramp over a duration of 1 s. Then, it is physically moved, such that it accelerates to 25 m/s along x-direction in 2 s, maintains this velocity for 10 s, and, then, decelerates to a stop in 1.2 s. In the ALE and MORALE simulations, a guiding velocity is supplied to the material of the structure (in the opposite direction of load movement corresponding to the Lagrangian simulation), such that an equivalent load movement (compared to the Lagrangian simulation) is achieved. In terms of discretization, uniform meshes with individual finite elements having dimensions 10 m $\times$ 10 m $\times$ 10 m, are adopted for Lagrangian, ALE and MORALE simulations. In other words, the finite element density is the same in all meshes. Further, care is taken to make comparison as fair as possible. Consider Figure \ref{fig:specimens00}(b), which shows the geometry of the specimen used in the Lagrangian simulation. The distance to the boundaries from the centre of the applied load is at minimum 80 m (along the x- and y-direction) over the length of the load path. The domain considered in the ALE and MORALE simulations is minimized, ensuring that this distance between the load centre and the boundary of 80 m (along the x- and y-direction), is maintained as shown in Figure \ref{fig:specimens00}(a). Further, the depth (along z-direction) of the specimens is also maintained at 80 m, in all three simulations. Accordingly, the Lagrangian simulation requires a discretization of the entire structure in the path of the moving load, therefore, 5760 finite elements are used, resulting in 15840 degrees of freedom. However, the ALE and MORALE simulations only require 2048 finite elements, resulting in 5400 degrees of freedom. The obtained displacement component $u_z$ under the load centre is plotted over time in Figure \ref{fig:dispvstime00}. Additionally, the times taken for the simulations to run through are listed in Table \ref{tab:timecomp} and plotted in Figure \ref{fig:runtimes00}. It is worth mentioning that in this study, all time dependent effects are only due to inertial effects because elastic materials are used. From the results shown in Figure \ref{fig:dispvstime00} and Table \ref{tab:timecomp}, it can be observed that the ALE simulations are significantly faster than the conventional Lagrangian simulation, while still achieving reasonable accuracy. The small discrepancies observed could possibly be attributed to errors during interpolation of the moving load in the Lagrangian simulation, and the time step and mesh sizes used. To apply the POD to the problem at hand, one full order ALE precomputation with the aforementioned simulation parameters is conducted, and all $\ell = 330$ solution vectors are used to construct the snapshot matrix. Truncation of the left singular vectors after the normalized singular values have decreased to $10^{-7}$ leads to a projection matrix with $m=73$ POD modes. Despite the relatively coarse mesh, the simulation time using the MORALE formulation reduces by a factor of $\sim$4.5 with respect to the full order ALE simulation and by a factor of $\sim$14 with respect to the full order Lagrangian simulation. As seen in Figure \ref{fig:dispvstime00}, the MORALE solution shows good agreement with the full order ALE solution. It should be noted that the time reduction factors are given for the online simulation phase. A beneficial speedup is therefore only achieved if the same reduced order model can be used to simulate altered problems with high accuracy, such that the methodology can be used, for e.g., uncertainty quantification. A further study is conducted, highlighting the prediction accuracy of the MORALE formulation. In this study, only snapshots of one full order ALE simulation with a maximum guiding velocity $w_{max} = 25$ m/s and load $p = 50$ Pa are used to create the ROM. Then, the reduced order model is used for simulations with different guiding velocities and loads. The results of this study are shown in Figure \ref{fig:POD_performance}. The ROM can successfully capture and represent the physical effects arising from varying the load and guiding velocity. If the guiding velocity is varied, the obtained displacements are also different due to inertial effects. Further, if the load is varied, the displacement also varies accordingly. Thus, the ROM demonstrates very good agreement with the reference simulations for all tested cases, and also for simulations where both the guiding velocity as well as the load are varied.

\begin{table}[h!]
\caption{Time comparison of the Lagrangian, POD-reduced Lagrangian, ALE, and POD-reduced ALE simulations.}
\label{tab:timecomp}
\begin{center}
\begin{tblr}{r|rrrr}
\hline
Time required for [s]& Lagrange & {POD-Lagrange \\ ($m=188$)} & ALE & {MORALE \\  ($m=73$)} \\
\hline
{Computation and assembly of \\ (reduced) material quantities} & {\\1101} & {\\1193} & {\\363} & {\\393} \\
{Solving (reduced) \\ system of equations} & {\\5385} & {\\103} & {\\1620} & {\\26} \\
{Computation of (reduced) \\ ALE quantities} & {\\-} & {\\-} & {\\26} & {\\45} \\
POD basis computation & - & 3 & - & 0.2 \\
\hline
Total & 6486 & 1299 & 2009 & 464.2 \\
Relative [-] & 1.0 & 0.2 & 0.31 & 0.07 \\
\hline
Peak RAM usage [MB] & 322 & 312 & 222 & 210
\end{tblr}
\end{center}
\end{table}

\subsection{Study 2: Complex example}\label{ssec:neovis}
\begin{figure}[h!]
	\centering
	\def\svgwidth{\textwidth}
	\begingroup%
	\makeatletter%
	\providecommand\color[2][]{%
		\errmessage{(Inkscape) Color is used for the text in Inkscape, but the package 'color.sty' is not loaded}%
		\renewcommand\color[2][]{}%
	}%
	\providecommand\transparent[1]{%
		\errmessage{(Inkscape) Transparency is used (non-zero) for the text in Inkscape, but the package 'transparent.sty' is not loaded}%
		\renewcommand\transparent[1]{}%
	}%
	\providecommand\rotatebox[2]{#2}%
	\newcommand*\fsize{\dimexpr\f@size pt\relax}%
	\newcommand*\lineheight[1]{\fontsize{\fsize}{#1\fsize}\selectfont}%
	\ifx\svgwidth\undefined%
	\setlength{\unitlength}{4012.28108083bp}%
	\ifx\svgscale\undefined%
	\relax%
	\else%
	\setlength{\unitlength}{\unitlength * \real{\svgscale}}%
	\fi%
	\else%
	\setlength{\unitlength}{\svgwidth}%
	\fi%
	\global\let\svgwidth\undefined%
	\global\let\svgscale\undefined%
	\makeatother%
	\begin{picture}(1,1.07839503)%
		\lineheight{1}%
		\setlength\tabcolsep{0pt}%
		\put(0,0){\includegraphics[width=\unitlength,page=1]{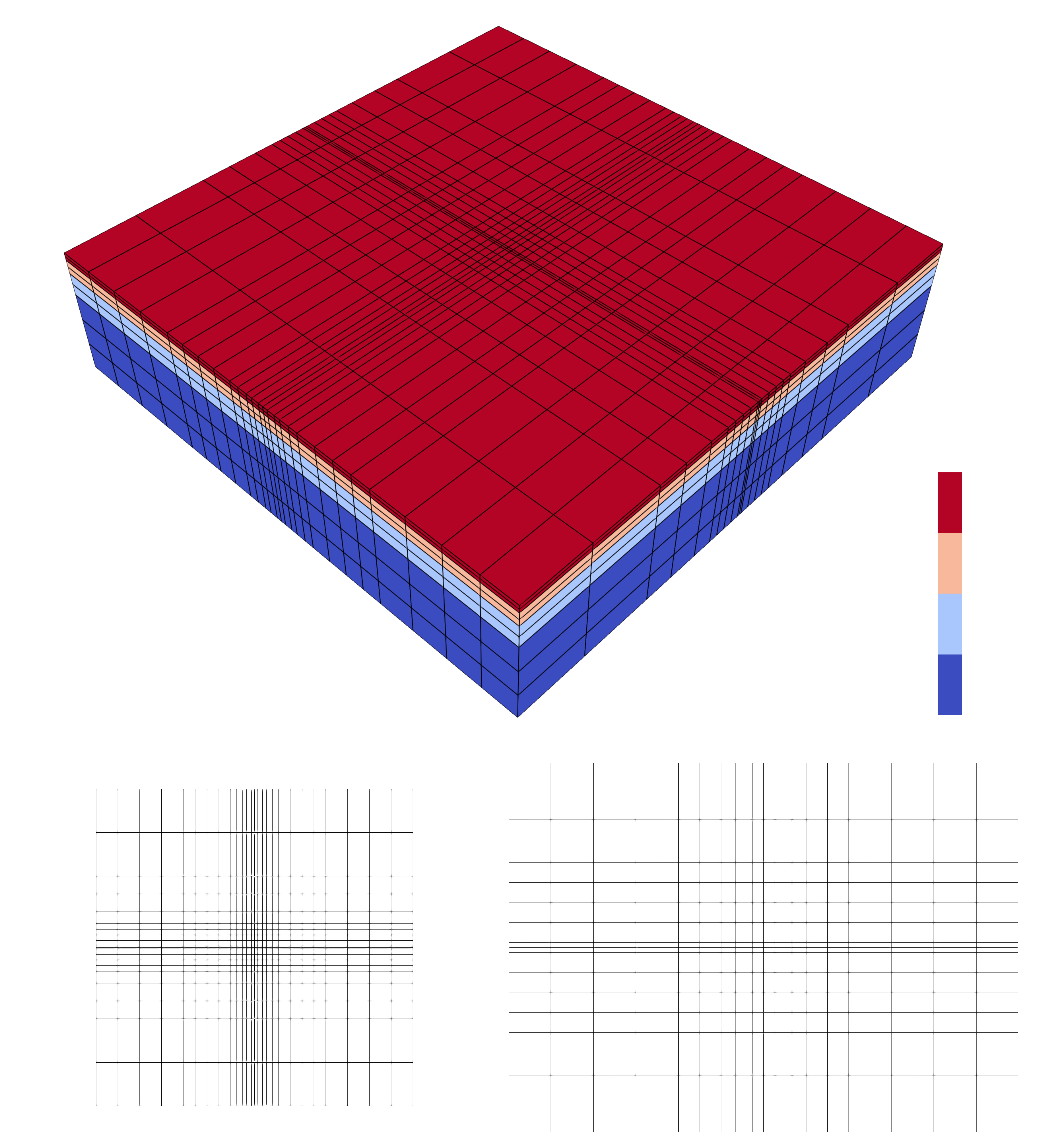}}%
		\put(0.03767731,1.04881852){\color[rgb]{0,0,0}\makebox(0,0)[lt]{\lineheight{1.25}\smash{\begin{tabular}[t]{l}(a)\end{tabular}}}}%
		\put(0.02070436,0.37610669){\color[rgb]{0,0,0}\makebox(0,0)[lt]{\lineheight{1.25}\smash{\begin{tabular}[t]{l}(b)\end{tabular}}}}%
		\put(0.44522663,0.37610669){\color[rgb]{0,0,0}\makebox(0,0)[lt]{\lineheight{1.25}\smash{\begin{tabular}[t]{l}(c)\end{tabular}}}}%
		\put(0,0){\includegraphics[width=\unitlength,page=2]{01_specimen_svg-tex.pdf}}%
		\put(0.1229291,0.48490175){\color[rgb]{0,0,0}\makebox(0,0)[lt]{\lineheight{1.25}\smash{\begin{tabular}[t]{l}X\end{tabular}}}}%
		\put(0.09381174,0.00571684){\color[rgb]{0,0,0}\makebox(0,0)[lt]{\lineheight{1.25}\smash{\begin{tabular}[t]{l}X\end{tabular}}}}%
		\put(0.12790443,0.56133573){\color[rgb]{0,0,0}\makebox(0,0)[lt]{\lineheight{1.25}\smash{\begin{tabular}[t]{l}Y\end{tabular}}}}%
		\put(0.01157727,0.08977504){\color[rgb]{0,0,0}\makebox(0,0)[lt]{\lineheight{1.25}\smash{\begin{tabular}[t]{l}Y\end{tabular}}}}%
		\put(0.06118709,0.58530431){\color[rgb]{0,0,0}\makebox(0,0)[lt]{\lineheight{1.25}\smash{\begin{tabular}[t]{l}Z\end{tabular}}}}%
		\put(0,0){\includegraphics[width=\unitlength,page=3]{01_specimen_svg-tex.pdf}}%
		\put(0.69260673,0.00155529){\color[rgb]{0,0,0}\makebox(0,0)[lt]{\lineheight{1.25}\smash{\begin{tabular}[t]{l}p = 0.818 MPa\end{tabular}}}}%
		\put(0,0){\includegraphics[width=\unitlength,page=4]{01_specimen_svg-tex.pdf}}%
		\put(0.70160123,0.38703741){\color[rgb]{0,0,0}\makebox(0,0)[lt]{\lineheight{1.25}\smash{\begin{tabular}[t]{l}160 mm\end{tabular}}}}%
		\put(0.22695863,0.36096354){\color[rgb]{0,0,0}\makebox(0,0)[lt]{\lineheight{1.25}\smash{\begin{tabular}[t]{l}8000 mm\end{tabular}}}}%
		\put(0.06141287,0.17521557){\color[rgb]{0,0,0}\rotatebox{90}{\makebox(0,0)[lt]{\lineheight{1.25}\smash{\begin{tabular}[t]{l}8000 mm\end{tabular}}}}}%
		\put(0.96942019,0.20960636){\color[rgb]{0,0,0}\makebox(0,0)[lt]{\lineheight{1.25}\smash{\begin{tabular}[t]{l}280 mm\end{tabular}}}}%
		\put(0.9712768,0.18558497){\color[rgb]{0,0,0}\makebox(0,0)[lt]{\lineheight{1.25}\smash{\begin{tabular}[t]{l}70 mm\end{tabular}}}}%
		\put(0.71583008,0.20974895){\color[rgb]{0,0,0}\makebox(0,0)[lt]{\lineheight{1.25}\smash{\begin{tabular}[t]{l}p\end{tabular}}}}%
		\put(0.71588993,0.1631719){\color[rgb]{0,0,0}\makebox(0,0)[lt]{\lineheight{1.25}\smash{\begin{tabular}[t]{l}p\end{tabular}}}}%
		\put(0.82752399,0.60461382){\color[rgb]{0,0,0}\makebox(0,0)[lt]{\lineheight{1.25}\smash{\begin{tabular}[t]{l}asphalt\end{tabular}}}}%
		\put(0.80412023,0.54637986){\color[rgb]{0,0,0}\makebox(0,0)[lt]{\lineheight{1.25}\smash{\begin{tabular}[t]{l}base course\end{tabular}}}}%
		\put(0.8213941,0.4895953){\color[rgb]{0,0,0}\makebox(0,0)[lt]{\lineheight{1.25}\smash{\begin{tabular}[t]{l}sub-base\end{tabular}}}}%
		\put(0.82959405,0.43281073){\color[rgb]{0,0,0}\makebox(0,0)[lt]{\lineheight{1.25}\smash{\begin{tabular}[t]{l}subsoil\end{tabular}}}}%
		\put(0.91358281,0.60457361){\color[rgb]{0,0,0}\makebox(0,0)[lt]{\lineheight{1.25}\smash{\begin{tabular}[t]{l}120 mm\end{tabular}}}}%
		\put(0.91385663,0.5473059){\color[rgb]{0,0,0}\makebox(0,0)[lt]{\lineheight{1.25}\smash{\begin{tabular}[t]{l}220 mm\end{tabular}}}}%
		\put(0.91383472,0.49003821){\color[rgb]{0,0,0}\makebox(0,0)[lt]{\lineheight{1.25}\smash{\begin{tabular}[t]{l}360 mm\end{tabular}}}}%
		\put(0.91358281,0.43277052){\color[rgb]{0,0,0}\makebox(0,0)[lt]{\lineheight{1.25}\smash{\begin{tabular}[t]{l}1300 mm\end{tabular}}}}%
	\end{picture}%
	\endgroup%
	\caption{(a) Perspective view of the non-uniform mesh showing layered structure, (b) top view of the mesh and (c) close-up view with loading details corresponding to simulations of Section \ref{ssec:neovis}.}
	\label{fig:specimen01}
\end{figure}

The simulations in this study use a nonlinear viscoelastic material model on a hyperelastic Neo-Hookean basis \cite{anantheswar2025treatment}. The material parameters used are listed in Table \ref{table:complex_params}. A non-uniform mesh with refinement in the loaded region is adopted, as shown in Figure \ref{fig:specimen01}.

\begin{table}[h!]
	\caption{Parameters chosen for the study in Section \ref{ssec:neovis}.}
	\begin{center}
		\begin{tblr}{colspec = {X[c,m]X[c,m]cX[c,m]c},
				cell{1}{2} = {c=2}{c}, 
				cell{1}{4} = {c=2}{c}, 
				cell{2}{2} = {c=1}{r},
				cell{2}{4} = {c=1}{r},
				cell{3}{2} = {c=1}{r},
                    cell{3}{4} = {c=1}{r},
				cell{4}{2} = {c=1}{r},
				cell{5}{2} = {c=1}{r},
				cell{5}{4} = {c=1}{r},
				cell{6}{2} = {c=1}{r},
                    cell{6}{4} = {c=1}{r},
				cell{7}{2} = {c=1}{r},
				cell{8}{2} = {c=1}{r},
				cell{8}{4} = {c=1}{r},
				cell{9}{2} = {c=1}{r},
                    cell{9}{4} = {c=1}{r},
				cell{10}{2} = {c=1}{r},
				cell{11}{2} = {c=1}{r},
				cell{11}{4} = {c=1}{r},
				cell{12}{2} = {c=1}{r},
                    cell{12}{4} = {c=1}{r},
				cell{13}{2} = {c=1}{r}
			}
			\hline
			layer &\vline elastic branch &  	&\vline viscous branch &	\\
			\hline
			asphalt & bulk modulus $\kappa$ [MPa] & 3.0$\times$10$^4$ & shear modulus $\mu_v$ [MPa] & 1.0$\times$10$^2$ \\
			& shear modulus $\mu$ [MPa] & 2.5$\times$10$^2$ & parameter $\eta_v$ [MPa s] & 1.0$\times$10$^3$ \\
			& density $\hat{\rho}$ [kg/m$^3$] & 2.3$\times$10$^{3}$ & &	\\
			\hline
			base course & bulk modulus $\kappa$ [MPa] & 4.5$\times$10$^2$ & shear modulus $\mu_v$ [MPa] & 5.0$\times$10$^1$ \\
			& shear modulus $\mu$ [MPa] & 1.5$\times$10$^2$ & parameter $\eta_v$ [MPa s] & 5.0$\times$10$^2$ \\
			& density $\hat{\rho}$ [kg/m$^3$] & 2.2$\times$10$^{3}$ & &	\\
			\hline
			sub-base & bulk modulus $\kappa$ [MPa] & 3.0$\times$10$^2$ & shear modulus $\mu_v$ [MPa] & 5.0$\times$10$^1$ \\
			& shear modulus $\mu$ [MPa] & 1.0$\times$10$^2$ & parameter $\eta_v$ [MPa s] & 3.333$\times$10$^2$ \\
			& density $\hat{\rho}$ [kg/m$^3$] & 2.1$\times$10$^{3}$ & &	\\
			\hline
			subsoil & bulk modulus $\kappa$ [MPa] & 1.0$\times$10$^2$ & shear modulus $\mu_v$ [MPa] & 1.0$\times$10$^1$ \\
			& shear modulus $\mu$ [MPa] & 3.75$\times$10$^1$ & parameter $\eta_v$ [MPa s] & 2.5$\times$10$^2$ \\
			& density $\hat{\rho}$ [kg/m$^3$] & 2.0$\times$10$^{3}$ & &	\\
                \hline
		\end{tblr}
	\end{center}
\label{table:complex_params}
\end{table}

\begin{figure}[h!]
	\centering
	\def\svgwidth{\textwidth}
	\begingroup%
	\makeatletter%
	\providecommand\color[2][]{%
		\errmessage{(Inkscape) Color is used for the text in Inkscape, but the package 'color.sty' is not loaded}%
		\renewcommand\color[2][]{}%
	}%
	\providecommand\transparent[1]{%
		\errmessage{(Inkscape) Transparency is used (non-zero) for the text in Inkscape, but the package 'transparent.sty' is not loaded}%
		\renewcommand\transparent[1]{}%
	}%
	\providecommand\rotatebox[2]{#2}%
	\newcommand*\fsize{\dimexpr\f@size pt\relax}%
	\newcommand*\lineheight[1]{\fontsize{\fsize}{#1\fsize}\selectfont}%
	\ifx\svgwidth\undefined%
	\setlength{\unitlength}{675bp}%
	\ifx\svgscale\undefined%
	\relax%
	\else%
	\setlength{\unitlength}{\unitlength * \real{\svgscale}}%
	\fi%
	\else%
	\setlength{\unitlength}{\svgwidth}%
	\fi%
	\global\let\svgwidth\undefined%
	\global\let\svgscale\undefined%
	\makeatother%
	\begin{picture}(1,0.66666667)%
		\lineheight{1}%
		\setlength\tabcolsep{0pt}%
		\put(0,0){\includegraphics[width=\unitlength,page=1]{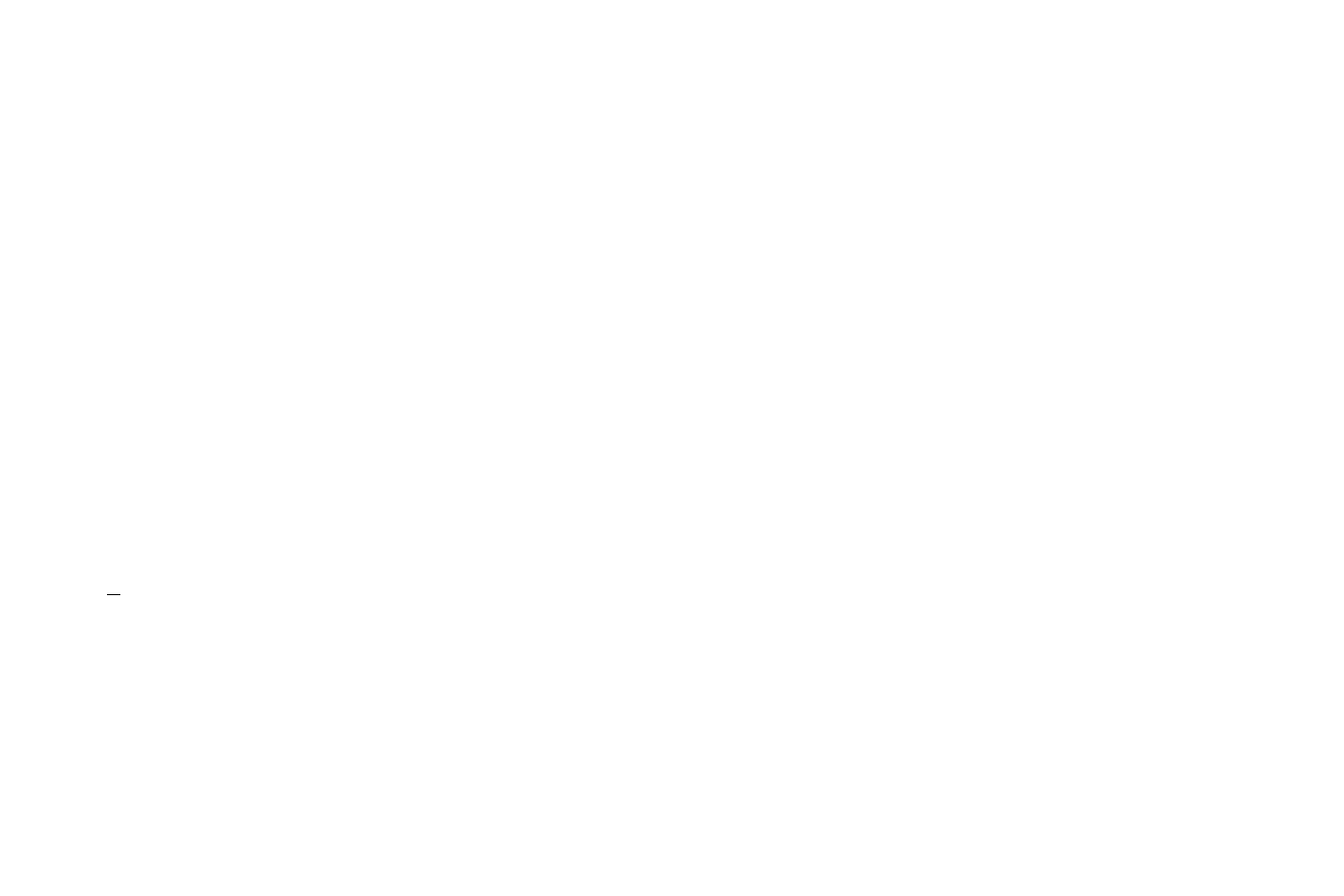}}%
		\put(0.07066667,0.21977778){\makebox(0,0)[rt]{\lineheight{1.25}\smash{\begin{tabular}[t]{r}-0.4\end{tabular}}}}%
		\put(0,0){\includegraphics[width=\unitlength,page=2]{01_disp_vs_time_svg-tex.pdf}}%
		\put(0.07066667,0.29322222){\makebox(0,0)[rt]{\lineheight{1.25}\smash{\begin{tabular}[t]{r}-0.3\end{tabular}}}}%
		\put(0,0){\includegraphics[width=\unitlength,page=3]{01_disp_vs_time_svg-tex.pdf}}%
		\put(0.07066667,0.36666667){\makebox(0,0)[rt]{\lineheight{1.25}\smash{\begin{tabular}[t]{r}-0.2\end{tabular}}}}%
		\put(0,0){\includegraphics[width=\unitlength,page=4]{01_disp_vs_time_svg-tex.pdf}}%
		\put(0.07066667,0.44022222){\makebox(0,0)[rt]{\lineheight{1.25}\smash{\begin{tabular}[t]{r}-0.1\end{tabular}}}}%
		\put(0,0){\includegraphics[width=\unitlength,page=5]{01_disp_vs_time_svg-tex.pdf}}%
		\put(0.07066667,0.51366667){\makebox(0,0)[rt]{\lineheight{1.25}\smash{\begin{tabular}[t]{r}0\end{tabular}}}}%
		\put(0,0){\includegraphics[width=\unitlength,page=6]{01_disp_vs_time_svg-tex.pdf}}%
		\put(0.07066667,0.58711111){\makebox(0,0)[rt]{\lineheight{1.25}\smash{\begin{tabular}[t]{r}0.1\end{tabular}}}}%
		\put(0,0){\includegraphics[width=\unitlength,page=7]{01_disp_vs_time_svg-tex.pdf}}%
		\put(0.07988889,0.14466667){\makebox(0,0)[t]{\lineheight{1.25}\smash{\begin{tabular}[t]{c}0\end{tabular}}}}%
		\put(0,0){\includegraphics[width=\unitlength,page=8]{01_disp_vs_time_svg-tex.pdf}}%
		\put(0.27911111,0.14466667){\makebox(0,0)[t]{\lineheight{1.25}\smash{\begin{tabular}[t]{c}5\end{tabular}}}}%
		\put(0,0){\includegraphics[width=\unitlength,page=9]{01_disp_vs_time_svg-tex.pdf}}%
		\put(0.47844444,0.14466667){\makebox(0,0)[t]{\lineheight{1.25}\smash{\begin{tabular}[t]{c}10\end{tabular}}}}%
		\put(0,0){\includegraphics[width=\unitlength,page=10]{01_disp_vs_time_svg-tex.pdf}}%
		\put(0.67766667,0.14466667){\makebox(0,0)[t]{\lineheight{1.25}\smash{\begin{tabular}[t]{c}15\end{tabular}}}}%
		\put(0,0){\includegraphics[width=\unitlength,page=11]{01_disp_vs_time_svg-tex.pdf}}%
		\put(0.877,0.14466667){\makebox(0,0)[t]{\lineheight{1.25}\smash{\begin{tabular}[t]{c}20\end{tabular}}}}%
		\put(0,0){\includegraphics[width=\unitlength,page=12]{01_disp_vs_time_svg-tex.pdf}}%
		\put(0.90611111,0.183){\makebox(0,0)[lt]{\lineheight{1.25}\smash{\begin{tabular}[t]{l}0\end{tabular}}}}%
		\put(0,0){\includegraphics[width=\unitlength,page=13]{01_disp_vs_time_svg-tex.pdf}}%
		\put(0.90611111,0.25655556){\makebox(0,0)[lt]{\lineheight{1.25}\smash{\begin{tabular}[t]{l}0.2\end{tabular}}}}%
		\put(0,0){\includegraphics[width=\unitlength,page=14]{01_disp_vs_time_svg-tex.pdf}}%
		\put(0.90611111,0.33){\makebox(0,0)[lt]{\lineheight{1.25}\smash{\begin{tabular}[t]{l}0.4\end{tabular}}}}%
		\put(0,0){\includegraphics[width=\unitlength,page=15]{01_disp_vs_time_svg-tex.pdf}}%
		\put(0.90611111,0.40344444){\makebox(0,0)[lt]{\lineheight{1.25}\smash{\begin{tabular}[t]{l}0.6\end{tabular}}}}%
		\put(0,0){\includegraphics[width=\unitlength,page=16]{01_disp_vs_time_svg-tex.pdf}}%
		\put(0.90611111,0.47688889){\makebox(0,0)[lt]{\lineheight{1.25}\smash{\begin{tabular}[t]{l}0.8\end{tabular}}}}%
		\put(0,0){\includegraphics[width=\unitlength,page=17]{01_disp_vs_time_svg-tex.pdf}}%
		\put(0.90611111,0.55033333){\makebox(0,0)[lt]{\lineheight{1.25}\smash{\begin{tabular}[t]{l}1\end{tabular}}}}%
		\put(0,0){\includegraphics[width=\unitlength,page=18]{01_disp_vs_time_svg-tex.pdf}}%
		\put(0.90611111,0.62388889){\makebox(0,0)[lt]{\lineheight{1.25}\smash{\begin{tabular}[t]{l}1.2\end{tabular}}}}%
		\put(0,0){\includegraphics[width=\unitlength,page=19]{01_disp_vs_time_svg-tex.pdf}}%
		\put(0.02111111,0.40777778){\rotatebox{90}{\makebox(0,0)[t]{\lineheight{1.25}\smash{\begin{tabular}[t]{c}Displacement $u_z$ [mm] under load\end{tabular}}}}}%
		\put(0.96655556,0.40777778){\rotatebox{90}{\makebox(0,0)[t]{\lineheight{1.25}\smash{\begin{tabular}[t]{c}Amplitude scaling factor [-]\end{tabular}}}}}%
		\put(0.48833333,0.10966667){\makebox(0,0)[t]{\lineheight{1.25}\smash{\begin{tabular}[t]{c}Time [s]\end{tabular}}}}%
		\put(0.34733333,0.08566667){\makebox(0,0)[rt]{\lineheight{1.25}\smash{\begin{tabular}[t]{r}Displacements:\end{tabular}}}}%
		\put(0.34733333,0.06566667){\makebox(0,0)[rt]{\lineheight{1.25}\smash{\begin{tabular}[t]{r}ALE\end{tabular}}}}%
		\put(0,0){\includegraphics[width=\unitlength,page=20]{01_disp_vs_time_svg-tex.pdf}}%
		\put(0.34733333,0.04566667){\makebox(0,0)[rt]{\lineheight{1.25}\smash{\begin{tabular}[t]{r}MORALE\end{tabular}}}}%
		\put(0,0){\includegraphics[width=\unitlength,page=21]{01_disp_vs_time_svg-tex.pdf}}%
		\put(0.68011111,0.08566667){\makebox(0,0)[rt]{\lineheight{1.25}\smash{\begin{tabular}[t]{r}Amplitude scaling factors:\end{tabular}}}}%
		\put(0.68011111,0.06566667){\makebox(0,0)[rt]{\lineheight{1.25}\smash{\begin{tabular}[t]{r}Pressure load\end{tabular}}}}%
		\put(0,0){\includegraphics[width=\unitlength,page=22]{01_disp_vs_time_svg-tex.pdf}}%
		\put(0.68011111,0.04566667){\makebox(0,0)[rt]{\lineheight{1.25}\smash{\begin{tabular}[t]{r}Load movement / Guiding velocity\end{tabular}}}}%
		\put(0,0){\includegraphics[width=\unitlength,page=23]{01_disp_vs_time_svg-tex.pdf}}%
	\end{picture}%
	\endgroup%
	\caption{Displacement $u_z$ under the centre of the load for the study in Section \ref{ssec:neovis}, plotted over time.}
	\label{fig:dispvstime01}
\end{figure}

\begin{figure}[h!]
	\centering
	\def\svgwidth{\textwidth}
	\begingroup%
	\makeatletter%
	\providecommand\color[2][]{%
		\errmessage{(Inkscape) Color is used for the text in Inkscape, but the package 'color.sty' is not loaded}%
		\renewcommand\color[2][]{}%
	}%
	\providecommand\transparent[1]{%
		\errmessage{(Inkscape) Transparency is used (non-zero) for the text in Inkscape, but the package 'transparent.sty' is not loaded}%
		\renewcommand\transparent[1]{}%
	}%
	\providecommand\rotatebox[2]{#2}%
	\newcommand*\fsize{\dimexpr\f@size pt\relax}%
	\newcommand*\lineheight[1]{\fontsize{\fsize}{#1\fsize}\selectfont}%
	\ifx\svgwidth\undefined%
	\setlength{\unitlength}{2042.04316868bp}%
	\ifx\svgscale\undefined%
	\relax%
	\else%
	\setlength{\unitlength}{\unitlength * \real{\svgscale}}%
	\fi%
	\else%
	\setlength{\unitlength}{\svgwidth}%
	\fi%
	\global\let\svgwidth\undefined%
	\global\let\svgscale\undefined%
	\makeatother%
	\begin{picture}(1,1.349167)%
		\lineheight{1}%
		\setlength\tabcolsep{0pt}%
		\put(0,0){\includegraphics[width=\unitlength,page=1]{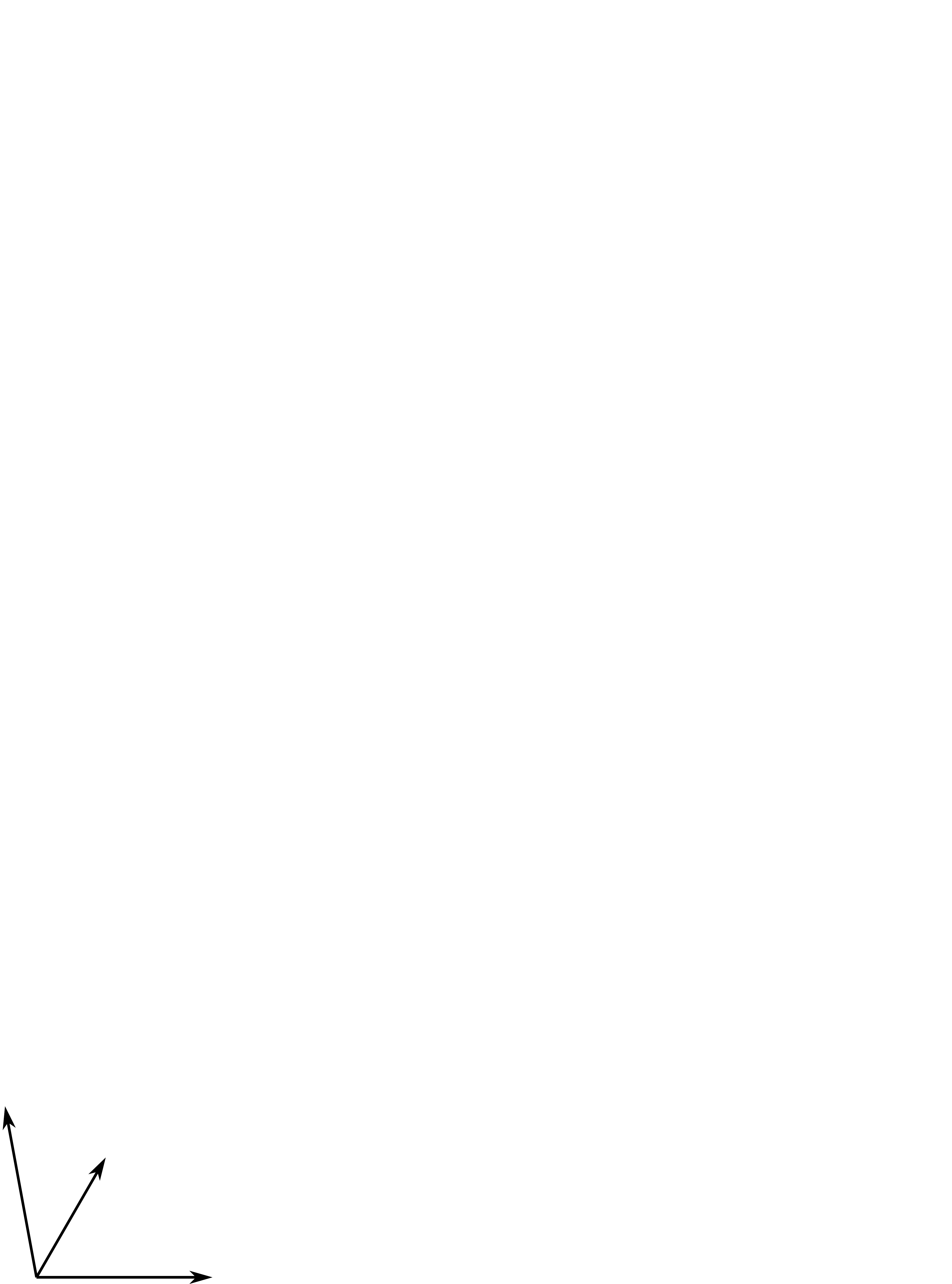}}%
		\put(0.23622916,0.00178145){\color[rgb]{0,0,0}\makebox(0,0)[lt]{\lineheight{1.25}\smash{\begin{tabular}[t]{l}X\end{tabular}}}}%
		\put(0.11626638,0.1395947){\color[rgb]{0,0,0}\makebox(0,0)[lt]{\lineheight{1.25}\smash{\begin{tabular}[t]{l}Y\end{tabular}}}}%
		\put(-0.00066,0.20085798){\color[rgb]{0,0,0}\makebox(0,0)[lt]{\lineheight{1.25}\smash{\begin{tabular}[t]{l}Z\end{tabular}}}}%
		\put(0.07027775,1.20557314){\makebox(0,0)[lt]{\lineheight{1.25}\smash{\begin{tabular}[t]{l}$t = 4.0$ s\end{tabular}}}}%
		\put(0.07027775,1.02391691){\makebox(0,0)[lt]{\lineheight{1.25}\smash{\begin{tabular}[t]{l}$t = 4.75$ s\end{tabular}}}}%
		\put(0.07027775,0.8422607){\makebox(0,0)[lt]{\lineheight{1.25}\smash{\begin{tabular}[t]{l}$t = 7.5$ s\end{tabular}}}}%
		\put(0.07027775,0.66060445){\makebox(0,0)[lt]{\lineheight{1.25}\smash{\begin{tabular}[t]{l}$t = 9.05$ s\end{tabular}}}}%
		\put(0.07027775,0.4789482){\makebox(0,0)[lt]{\lineheight{1.25}\smash{\begin{tabular}[t]{l}$t = 9.15$ s\end{tabular}}}}%
		\put(0.07027897,0.29728745){\makebox(0,0)[lt]{\lineheight{1.25}\smash{\begin{tabular}[t]{l}$t = 18.0$ s\end{tabular}}}}%
		\put(0.53672189,0.14681305){\makebox(0,0)[lt]{\lineheight{1.25}\smash{\begin{tabular}[t]{l}$e_{zz}$ [-]\end{tabular}}}}%
		\put(0.24898068,0.10481032){\makebox(0,0)[lt]{\lineheight{1.25}\smash{\begin{tabular}[t]{l}2.5e-5\end{tabular}}}}%
		\put(0.1553421,0.10481032){\makebox(0,0)[lt]{\lineheight{1.25}\smash{\begin{tabular}[t]{l}9.5e-5\end{tabular}}}}%
		\put(0.32773418,0.10481032){\makebox(0,0)[lt]{\lineheight{1.25}\smash{\begin{tabular}[t]{l}-4.5e-5\end{tabular}}}}%
		\put(0.40887897,0.10481032){\makebox(0,0)[lt]{\lineheight{1.25}\smash{\begin{tabular}[t]{l}-1.2e-4\end{tabular}}}}%
		\put(0.49026862,0.10481032){\makebox(0,0)[lt]{\lineheight{1.25}\smash{\begin{tabular}[t]{l}-1.8e-4\end{tabular}}}}%
		\put(0.57165809,0.10481032){\makebox(0,0)[lt]{\lineheight{1.25}\smash{\begin{tabular}[t]{l}-2.6e-4\end{tabular}}}}%
		\put(0.65304757,0.10481032){\makebox(0,0)[lt]{\lineheight{1.25}\smash{\begin{tabular}[t]{l}-3.2e-4\end{tabular}}}}%
		\put(0.73443722,0.10481032){\makebox(0,0)[lt]{\lineheight{1.25}\smash{\begin{tabular}[t]{l}-4.0e-4\end{tabular}}}}%
		\put(0.81582686,0.10481032){\makebox(0,0)[lt]{\lineheight{1.25}\smash{\begin{tabular}[t]{l}-4.7e-4\end{tabular}}}}%
		\put(0.88056186,0.10481032){\makebox(0,0)[lt]{\lineheight{1.25}\smash{\begin{tabular}[t]{l}-5.0e-4\end{tabular}}}}%
		\put(0,0){\includegraphics[width=\unitlength,page=2]{01_strain_contours_svg-tex.pdf}}%
		\put(0.3120964,1.32774712){\makebox(0,0)[lt]{\lineheight{1.25}\smash{\begin{tabular}[t]{l}ALE simulation\end{tabular}}}}%
		\put(0.7020077,1.32776147){\makebox(0,0)[lt]{\lineheight{1.25}\smash{\begin{tabular}[t]{l}MORALE simulation\end{tabular}}}}%
		\put(0,0){\includegraphics[width=\unitlength,page=3]{01_strain_contours_svg-tex.pdf}}%
	\end{picture}%
	\endgroup%
	\caption{Euler-Almansi strain $e_{zz}$ contours for the ALE and the MORALE simulations in Section \ref{ssec:neovis}, plotted on one half (cut along x-direction) of the specimen, at various time points during the simulations.}
	\label{fig:straincont01}
\end{figure}

This non-uniform mesh enables the application of a realistic wheel load from a dual truck tire. With the conventional Lagrangian framework, a refined mesh would need to be used in the entire path of the moving load, and this would be very inefficient. Therefore, only ALE and MORALE simulations are compared in this study. The boundary conditions are such that all surfaces of the structure, except for the top surface (where the load is applied), are fully restrained, resulting in 12825 degrees of freedom. Amplitude scaling factors for the load application and guiding velocity are plotted in Figure \ref{fig:dispvstime01}. The loading is first applied as a linear ramp over a period of 2 s. Then, the guiding velocity is applied, such that the load appears to accelerate to a velocity of 100 km/h in 3 s, starting at time $t = 4$ s. The amplitude of the load is then varied by $\pm$20\% in the form of a triangular wave, with a frequency of 5 Hz between $t = 8$ s and $t = 12$ s, while maintaining the velocity of 100 km/h. Following this, the load amplitude is maintained constant for the rest of the simulation. Finally, a deceleration is applied from 100 km/h at $t = 14$ s to a stop at $t = 19$ s. The resulting displacements at the central point on the top surface of the mesh are plotted over time in Figure \ref{fig:dispvstime01} for the ALE and MORALE simulations. Contour plots of the Euler-Almansi strain component $e_{zz}$ are depicted at various stages of the full order ALE simulation and the MORALE simulation on one half (cut along the x-direction) of the specimen in Figure \ref{fig:straincont01}.

\begin{table}[h!]
\caption{Time comparison of the ALE, and POD-reduced ALE simulations for the complex example.}
\label{tab:timecomp2}
\begin{center}
\begin{tblr}{r|rr}
\hline
Time required for [s] & ALE & {MORALE \\  ($m=168$)} \\
\hline
{Computation and assembly of (reduced) quantities} & {3635} & {3761} \\
{Solving (reduced) system of equations} & {22123} & {104} \\
POD basis computation & - & 8 \\
\hline
Total & 25758 & 3873 \\
Relative [-] & 1.0 & 0.15 \\
\hline
Peak RAM usage [MB] & 580 & 556
\end{tblr}
\end{center}
\end{table}

\begin{figure}[h!]
	\centering
	\def\svgwidth{\textwidth}
	\begingroup%
	\makeatletter%
	\providecommand\color[2][]{%
		\errmessage{(Inkscape) Color is used for the text in Inkscape, but the package 'color.sty' is not loaded}%
		\renewcommand\color[2][]{}%
	}%
	\providecommand\transparent[1]{%
		\errmessage{(Inkscape) Transparency is used (non-zero) for the text in Inkscape, but the package 'transparent.sty' is not loaded}%
		\renewcommand\transparent[1]{}%
	}%
	\providecommand\rotatebox[2]{#2}%
	\newcommand*\fsize{\dimexpr\f@size pt\relax}%
	\newcommand*\lineheight[1]{\fontsize{\fsize}{#1\fsize}\selectfont}%
	\ifx\svgwidth\undefined%
	\setlength{\unitlength}{460.79998779bp}%
	\ifx\svgscale\undefined%
	\relax%
	\else%
	\setlength{\unitlength}{\unitlength * \real{\svgscale}}%
	\fi%
	\else%
	\setlength{\unitlength}{\svgwidth}%
	\fi%
	\global\let\svgwidth\undefined%
	\global\let\svgscale\undefined%
	\makeatother%
	\begin{picture}(1,0.75000003)%
		\lineheight{1}%
		\setlength\tabcolsep{0pt}%
		\put(0,0){\includegraphics[width=\unitlength,page=1]{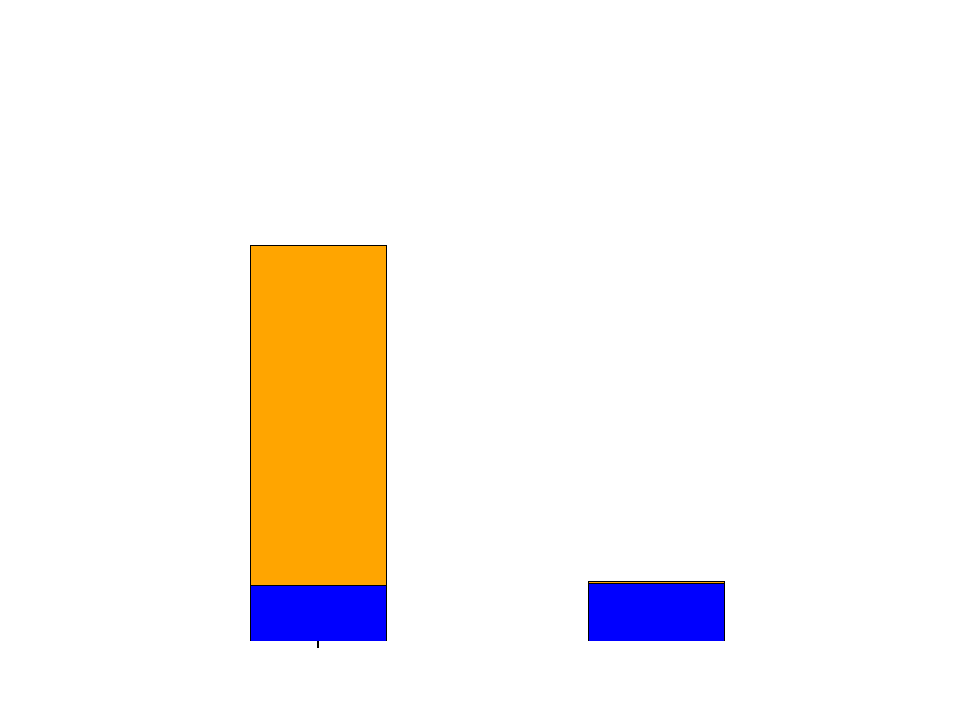}}%
		\put(0.33159328,0.05081939){\makebox(0,0)[t]{\lineheight{1.25}\smash{\begin{tabular}[t]{c}ALE\end{tabular}}}}%
		\put(0,0){\includegraphics[width=\unitlength,page=2]{2_timings_svg-tex.pdf}}%
		\put(0.65132805,0.05081939){\makebox(0,0)[lt]{\lineheight{1.25}\smash{\begin{tabular}[t]{l}MORALE\end{tabular}}}}%
		\put(0.65716543,0.02740019){\makebox(0,0)[lt]{\lineheight{1.25}\smash{\begin{tabular}[t]{l}($m$ = 168)\end{tabular}}}}%
		\put(0,0){\includegraphics[width=\unitlength,page=3]{2_timings_svg-tex.pdf}}%
		\put(0.10980903,0.0742552){\makebox(0,0)[rt]{\lineheight{1.25}\smash{\begin{tabular}[t]{r}0.0\end{tabular}}}}%
		\put(0,0){\includegraphics[width=\unitlength,page=4]{2_timings_svg-tex.pdf}}%
		\put(0.10980903,0.15675519){\makebox(0,0)[rt]{\lineheight{1.25}\smash{\begin{tabular}[t]{r}0.2\end{tabular}}}}%
		\put(0,0){\includegraphics[width=\unitlength,page=5]{2_timings_svg-tex.pdf}}%
		\put(0.10980903,0.2392552){\makebox(0,0)[rt]{\lineheight{1.25}\smash{\begin{tabular}[t]{r}0.4\end{tabular}}}}%
		\put(0,0){\includegraphics[width=\unitlength,page=6]{2_timings_svg-tex.pdf}}%
		\put(0.10980903,0.32175519){\makebox(0,0)[rt]{\lineheight{1.25}\smash{\begin{tabular}[t]{r}0.6\end{tabular}}}}%
		\put(0,0){\includegraphics[width=\unitlength,page=7]{2_timings_svg-tex.pdf}}%
		\put(0.10980903,0.4042552){\makebox(0,0)[rt]{\lineheight{1.25}\smash{\begin{tabular}[t]{r}0.8\end{tabular}}}}%
		\put(0,0){\includegraphics[width=\unitlength,page=8]{2_timings_svg-tex.pdf}}%
		\put(0.10980903,0.48675519){\makebox(0,0)[rt]{\lineheight{1.25}\smash{\begin{tabular}[t]{r}1.0\end{tabular}}}}%
		\put(0,0){\includegraphics[width=\unitlength,page=9]{2_timings_svg-tex.pdf}}%
		\put(0.10980903,0.5692552){\makebox(0,0)[rt]{\lineheight{1.25}\smash{\begin{tabular}[t]{r}1.2\end{tabular}}}}%
		\put(0,0){\includegraphics[width=\unitlength,page=10]{2_timings_svg-tex.pdf}}%
		\put(0.10980903,0.6517552){\makebox(0,0)[rt]{\lineheight{1.25}\smash{\begin{tabular}[t]{r}1.4\end{tabular}}}}%
		\put(0.06210327,0.37125002){\rotatebox{90}{\makebox(0,0)[t]{\lineheight{1.25}\smash{\begin{tabular}[t]{c}Relative simulation time [-]\end{tabular}}}}}%
		\put(0,0){\includegraphics[width=\unitlength,page=11]{2_timings_svg-tex.pdf}}%
		\put(0.34839141,0.62397911){\makebox(0,0)[lt]{\lineheight{1.25}\smash{\begin{tabular}[t]{l}Computing and assembly of (reduced) quantities\end{tabular}}}}%
		\put(0,0){\includegraphics[width=\unitlength,page=12]{2_timings_svg-tex.pdf}}%
		\put(0.34839141,0.59212554){\makebox(0,0)[lt]{\lineheight{1.25}\smash{\begin{tabular}[t]{l}Solving (reduced) system of equations\end{tabular}}}}%
		\put(0,0){\includegraphics[width=\unitlength,page=13]{2_timings_svg-tex.pdf}}%
		\put(0.34839141,0.56027197){\makebox(0,0)[lt]{\lineheight{1.25}\smash{\begin{tabular}[t]{l}POD basis computation\end{tabular}}}}%
	\end{picture}%
	\endgroup%
	\caption{Comparison of times required for simulations in Section \ref{ssec:neovis}.}
	\label{fig:runtimes01}
\end{figure}

From the results, it is apparent that the dynamic ALE formulation can satisfactorily capture the transient response of a pavement structure subjected to a moving wheel load. Further, the total distance traversed by the wheel load during this simulation is about 305.56 m. The benefit of applying the dynamic ALE formulation is that this entire length does not need be considered (contrary to a conventional Lagrangian simulation), instead, only the relevant region (8 m along both x- and y-directions and 2 m along the z-direction) around the load needs to be considered. This enables the simulation of complex multilayered structures like pavements within a reasonable time frame. It is worth mentioning that a time step size convergence study has been carried out on this example (omitted here for the sake of brevity). From this study, it has been observed that the obtained pavement response (see Figure \ref{fig:dispvstime01}) shows no significant dependence on the time step size. However, this is true only for the particular set of parameters and mesh chosen in this study. In general, care should be taken to ensure that time step size is sufficiently small, such that during the advection process, material points can be subjected to the load (without excessive skipping), and the internal variables can evolve. This improves the accuracy of the advection process for internal variables.

Further, applying the POD with $m=168$ modes, the simulation time is reduced by a factor of $\sim$6.7 as shown in Table \ref{tab:timecomp2} and Figure \ref{fig:runtimes01}. Overall, a good agreement between the result of the full order ALE simulation and the MORALE simulation can be observed, as shown in Figure \ref{fig:dispvstime01}. Only slight prediction inaccuracies can be observed, when the deceleration starts at $t=14$ s. Further, the strain contours plotted in Figure \ref{fig:straincont01} show satisfactory agreement between the ALE and MORALE simulations, however, the MORALE simulation does show some noise in the lower strain range. This could be caused by the effects of resonance in the MORALE simulation (because a lower number of modes is considered). However, further investigations are necessary.

\section{Conclusion}\label{sec:conc}
In this contribution, a novel Model Order Reduced Arbitrary Lagrangian Eulerian (MORALE) formulation for the transient simulation of structures like pavements, which are subjected to moving loads, has been put forth. It has been demonstrated through numerical studies that the formulation is faster than the Arbitrary Lagrangian Eulerian (ALE) approach alone, faster than the Lagrangian approach with Proper Orthogonal Decomposition (POD), and also significantly faster than the conventional Lagrangian approach. Further, the MORALE approach has been demonstrated to work well, even with inelastic materials and non-uniform meshes. Despite provision of a limited amount of training data, the MORALE approach achieves high prediction accuracy, reducing the computational effort significantly. The framework could therefore be used to boost the efficiency of tasks such as sensitivity analysis or uncertainty quantification. To increase the accuracy of the MORALE approach further, it could be beneficial to apply clustering techniques to the snapshots and create local ROMs, specifically tailored to e.g., the acceleration, cruising or deceleration phases. \\
Future research should be directed towards the advancement of the presented model order reduction technique through hyper-reduction strategies like the Discrete Empirical Interpolation Method (DEIM) or the Energy-Conserving Sampling and Weighting method (ECSW). Currently, these advanced MOR techniques are not directly applicable with ALE approaches because of the requirements of the internal variable update procedures. Further developments along this direction will likely further increase the efficiency and speed of computations, which is vital in digital twins, because the twinning process necessitates fast (almost real-time) simulations. Further developments could include the extension of the MORALE framework to utilize advanced multiphysics based material models for the various layers of the pavement structure, to enable improved prediction capabilities. Analysis of the performance of the developed MORALE formulation, when applied to more complex geometries with advanced material models (for example without the assumption of longitudinal homogeneity), could lead to further improvements. Moreover, extension of the formulation, such that it can be applied to rolling tires, is another avenue worth exploring.

\section*{Acknowledgements}\label{sec:ack}
The present work has been developed under subprojects A01 and B05 of the research project No. 453596084 (SFB/TRR 339), which has been granted by the German Research Foundation (Deutsche Forschungsgemeinschaft). This financial support is gratefully acknowledged.

\section*{Data Availability Statement} The data that support the findings of this study are available from the authors upon reasonable request.

\bibliography{bibliography}

@book{wriggers2008nonlinear,
	title={Nonlinear finite element methods},
	author={Wriggers, Peter},
	year={(2008)},
	address={Berlin, Germany},
	publisher={Springer Science \& Business Media}
}

@book{holzapfel2002nonlinear,
	title={Nonlinear solid mechanics: a continuum approach for engineering},
	author={Holzapfel, GA},
	year={(2000)},
	address={Chichester, England},
	publisher={Wiley}
}

@article{taylor_feap_2020,
	title = {{FEAP} -- {F}inite {E}lement {A}nalysis {P}rogram, url: http://projects.ce.berkeley.edu/feap, University of California, Berkeley},
	url = {http://projects.ce.berkeley.edu/feap},
	publisher = {University of California, Berkeley},
	author = {Taylor, RL},
	year = {(2020)}
}

@article{zopf2015continuum,
	title={A continuum mechanical approach to model asphalt},
	author={Zopf, C and Garcia, MA and Kaliske, M},
	journal={International Journal of Pavement Engineering},
	volume={16},
	pages={105--124},
	year={(2015)},
	publisher={Taylor \& Francis}
}

@article{kaliske2023digitaler,
	title={{D}igitaler {Z}willing {S}tra{\ss}e},
	author={Kaliske, Michael and Wollny, Ines},
	journal={Bautechnik},
	volume={100},
	pages={374--382},
	year={(2023)},
	publisher={Wiley Online Library}
}

@article{kaliske2022weg,
	title={{W}elchen {W}eg nimmt die „{S}tra{\ss}e der {Z}ukunft “? {D}igitalisierung der {S}tra{\ss}e im {S}onderforschungsbereich/{T}ransregio 339 „{D}igitaler {Z}willing {S}tra{\ss}e “},
	author={Kaliske, M and Oeser, M and Wollny, I and Behnke, R},
	journal={Bauingenieur},
	volume={97},
	pages={29--37},
	year={(2022)},
	publisher={ }
}

@article{venkatasubban1995new,
	title={A new finite element formulation for {ALE (arbitrary Lagrangian Eulerian)} compressible fluid mechanics},
	author={Venkatasubban, Chittur},
	journal={International Journal of Engineering Science},
	volume={33},
	pages={1743--1762},
	year={(1995)},
	publisher={Elsevier}
}

@article{souli2000ale,
	title={{ALE} formulation for fluid--structure interaction problems},
	author={Souli, Mohamed and Ouahsine, A and Lewin, L},
	journal={Computer Methods in Applied Mechanics and Engineering},
	volume={190},
	pages={659--675},
	year={(2000)},
	publisher={Elsevier}
}

@article{codina2009fixed,
	title={The fixed-mesh {ALE} approach for the numerical approximation of flows in moving domains},
	author={Codina, Ramon and Houzeaux, Guillaume and Coppola-Owen, Herbert and Baiges, Joan},
	journal={Journal of Computational Physics},
	volume={228},
	pages={1591--1611},
	year={(2009)},
	publisher={Elsevier}
}

@article{basting2017extended,
	title={Extended {ALE} method for fluid--structure interaction problems with large structural displacements},
	author={Basting, Steffen and Quaini, Annalisa and {\v{C}}ani{\'c}, Sun{\v{c}}ica and Glowinski, Roland},
	journal={Journal of Computational Physics},
	volume={331},
	pages={312--336},
	year={(2017)},
	publisher={Elsevier}
}

@article{benson1989efficient,
	title={An efficient, accurate, simple {ALE} method for nonlinear finite element programs},
	author={Benson, David},
	journal={Computer Methods in Applied Mechanics and Engineering},
	volume={72},
	pages={305--350},
	year={(1989)},
	publisher={Elsevier}
}

@article{bayoumi2004complete,
	title={A complete finite element treatment for the fully coupled implicit {ALE} formulation},
	author={Bayoumi, HN and Gadala, MohamedS},
	journal={Computational Mechanics},
	volume={33},
	pages={435--452},
	year={(2004)},
	publisher={Springer}
}

@inbook{donea2004arbitrary,
	title={{A}rbitrary {L}agrangian--{E}ulerian {M}ethods},
	author={Donea, Jean and Huerta, Antonio and Ponthot, J and Rodr{\'\i}guez-Ferran, Antonio},
	address={Ch. 14; Encyclopedia of Computational Mechanics},
	publisher={Wiley Online Library. (2004)}
}

@article{nazem2009arbitrary,
	title={{A}rbitrary {L}agrangian--{E}ulerian method for dynamic analysis of geotechnical problems},
	author={Nazem, Majidreza and Carter, John and Airey, David},
	journal={Computers and Geotechnics},
	volume={36},
	pages={549--557},
	year={(2009)},
	publisher={Elsevier}
}

@article{zreid2021ale,
	title={{ALE} formulation for thermomechanical inelastic material models applied to tire forming and curing simulations},
	author={Zreid, Imadeddin and Behnke, Ronny and Kaliske, Michael},
	journal={Computational Mechanics},
	volume={67},
	pages={1543--1557},
	year={(2021)},
	publisher={Springer}
}

@article{berger2022arbitrary,
	title={An arbitrary {L}agrangian {E}ulerian formulation for tire production simulation},
	author={Berger, Thomas and Kaliske, Michael},
	journal={Finite Elements in Analysis and Design},
	volume={204},
	pages={103742},
	year={(2022)},
	publisher={Elsevier}
}

@article{nackenhorst2004ale,
	title={The {ALE}-formulation of bodies in rolling contact: {T}heoretical foundations and finite element approach},
	author={Nackenhorst, U},
	journal={Computer Methods in Applied Mechanics and Engineering},
	volume={193},
	pages={4299--4322},
	year={(2004)},
	publisher={Elsevier}
}

@article{ziefle2008numerical,
	title={Numerical techniques for rolling rubber wheels: treatment of inelastic material properties and frictional contact},
	author={Ziefle, Matthias and Nackenhorst, Udo},
	journal={Computational Mechanics},
	volume={42},
	pages={337--356},
	year={(2008)},
	publisher={Springer}
}

@article{wollny2013numerical,
	title={Numerical simulation of pavement structures with inelastic material behaviour under rolling tyres based on an arbitrary {L}agrangian {E}ulerian ({ALE}) formulation},
	author={Wollny, Ines and Kaliske, Michael},
	journal={Road Materials and Pavement Design},
	volume={14},
	pages={71--89},
	year={(2013)},
	publisher={Taylor \& Francis}
}

@article{wollny2016numerical,
	title={Numerical modeling of inelastic structures at loading of steady state rolling},
	author={Wollny, Ines and Hartung, Felix and Kaliske, Michael},
	journal={Computational Mechanics},
	volume={57},
	pages={867--886},
	year={(2016)},
	publisher={Springer}
}

@article{anantheswar2023dynamic,
	title={A dynamic {ALE} formulation for structures under moving loads},
	author={Anantheswar, Atul and Wollny, Ines and Kaliske, Michael},
	journal={Computational Mechanics},
	volume={73},
	pages={139--157},
	year={(2024)},
	publisher={Springer}
}

@article{anantheswar2023transient,
	title={Transient response of pavement structures under moving wheel loads using the {ALE} methodology},
	author={Anantheswar, Atul and Wollny, Ines and Kaliske, Michael},
	journal={Proceedings in Applied Mathematics and Mechanics},
	volume={23},
	pages={e202300249},
	year={(2023)},
	publisher={Wiley Online Library}
}

@article{anantheswar2025treatment,
	title={Treatment of inelastic material models within a dynamic {ALE} formulation for structures subjected to moving loads},
	author={Anantheswar, Atul and Wollny, Ines and Kaliske, Michael},
	journal={International Journal for Numerical Methods in Engineering},
	volume={126},
	pages={e7599},
	year={(2025)},
	publisher={Wiley Online Library}
}

@article{wang2024nonsmooth,
  title={Nonsmooth model order reduction for transient tire--road dynamics of frictional contact with {ALE} formulations},
  author={Wang, Kun and Luo, Kai and Tian, Qiang},
  journal={Nonlinear Dynamics},
  volume={112},
  pages={18847--18868},
  year={(2024)},
  publisher={Springer}
}

@article{de2019development,
  title={Development and validation of a fully predictive high-fidelity simulation approach for predicting coarse road dynamic tire/road rolling contact forces},
  author={De Gregoriis, Daniel and Naets, Frank and Kindt, Peter and Desmet, Wim},
  journal={Journal of Sound and Vibration},
  volume={452},
  pages={147--168},
  year={(2019)},
  publisher={Elsevier}
}

@article{de2019application,
  title={Application of a priori hyper-reduction to the nonlinear dynamic finite element simulation of a rolling car tire},
  author={De Gregoriis, Daniel and Naets, Frank and Kindt, Peter and Desmet, Wim},
  journal={Journal of Computational and Nonlinear Dynamics},
  volume={14},
  pages={111009},
  year={(2019)},
  publisher={American Society of Mechanical Engineers}
}

@article{am3p_paper_jannick,
    author = {Kehls, Jannick and Anantheswar, Atul and Brepols, Tim and Wollny, Ines and Kaliske, Michael and Reese, Stefanie},
    title = {Towards real-time structural simulations for the digital twin of the road},
    journal = {Proceedings in Advances in Materials and Pavement Performance Prediction},
    year = {(2025)},
    volume = {4},
    pages = {596--599},
    url = {https://repositum.tuwien.at/handle/20.500.12708/214864}
}

@article{am3p_paper_atul,
    author = {Anantheswar, Atul and Wollny, Ines and Kaliske, Michael},
    title = {Simulating multilayered inelastic pavements by a dynamic {ALE} formulation},
    journal = {Proceedings in Advances in Materials and Pavement Performance Prediction},
    year = {(2025)},
    volume = {4},
    pages = {588--591},
    url = {https://repositum.tuwien.at/handle/20.500.12708/214864}
}

@article{newmark1959method,
	title={A method of computation for structural dynamics},
	author={Newmark, NathanM},
	journal={Journal of the Engineering Mechanics Division},
	volume={85},
	pages={67--94},
	year={(1959)},
	publisher={American Society of Civil Engineers}
}

@article{radermacher_comparison_2013,
	title = {A comparison of projection-based model reduction concepts in the context of nonlinear biomechanics},
	volume = {83},
	issn = {0939-1533, 1432-0681},
	url = {http://link.springer.com/10.1007/s00419-013-0742-9},
	language = {en},
	urldate = {2023-05-08},
	journal = {Archive of Applied Mechanics},
	author = {Radermacher, Annika and Reese, Stefanie},
	month = aug,
	year = {(2013)},
	pages = {1193--1213}
}

@book{hesthaven_certified_2016,
	title = {Certified {Reduced} {Basis} {Methods} for {Parametrized} {Partial} {Differential} {Equations}},
	isbn = {978-3-319-22469-5 978-3-319-22470-1},
	url = {https://link.springer.com/10.1007/978-3-319-22470-1},
	language = {en},
	urldate = {2023-05-08},
	author = {Hesthaven, S and Rozza, Gianluigi and Stamm, Benjamin},
	address = {Cham},
	publisher = {Springer},
	year = {(2016)},
	doi = {10.1007/978-3-319-22470-1}
}

@article{benner_survey_2015,
	title = {A {Survey} of {Projection}-{Based} {Model} {Reduction} {Methods} for {Parametric} {Dynamical} {Systems}},
	volume = {57},
	issn = {0036-1445, 1095-7200},
	url = {http://epubs.siam.org/doi/10.1137/130932715},
	language = {en},
	urldate = {2023-05-08},
	journal = {SIAM Review},
	author = {Benner, Peter and Gugercin, Serkan and Willcox, Karen},
	month = jan,
	year = {(2015)},
	pages = {483--531}
}

@article{spiess_reduction_2005,
	title = {Reduction {Methods} for {FE} {Analysis} in {Nonlinear} {Structural} {Dynamics}},
	volume = {5},
	copyright = {Copyright © 2005 WILEY-VCH Verlag GmbH \& Co. KGaA, Weinheim},
	issn = {1617-7061},
	url = {https://onlinelibrary.wiley.com/doi/abs/10.1002/pamm.200510048},
	language = {en},
	urldate = {2024-08-06},
	journal = {Proceedings in Applied Mathematics and Mechanics},
	author = {Spiess, Holger and Wriggers, Peter},
	year = {(2005)},
	pages = {135--136}
}

@article{gugercin_survey_2004,
	title = {A {Survey} of {Model} {Reduction} by {Balanced} {Truncation} and {Some} {New} {Results}},
	volume = {77},
	issn = {0020-7179},
	url = {https://doi.org/10.1080/00207170410001713448},
	urldate = {2025-04-02},
	journal = {International Journal of Control},
	author = {Gugercin, Serkan and Antoulas, C},
	month = may,
	year = {(2004)},
        publisher = {Taylor \& Francis},
	pages = {748--766}
}

@article{ritzert2025component,
  title={Component-Based Model-Order Reduction With Mortar Tied Contact for Nonlinear Quasi-Static Mechanical Problems},
  author={Ritzert, Stephan and Kehls, Jannick and Reese, Stefanie and Brepols, Tim},
  journal={International Journal for Numerical Methods in Engineering},
  volume={126},
  pages={e70041},
  year={(2025)},
  publisher={Wiley Online Library}
}

@article{zhang2025multi,
  title={A multi-field decomposed model order reduction approach for thermo-mechanically coupled gradient-extended damage simulations},
  author={Zhang, Qinghua and Ritzert, Stephan and Zhang, Jian and Kehls, Jannick and Reese, Stefanie and Brepols, Tim},
  journal={Computer Methods in Applied Mechanics and Engineering},
  volume={434},
  pages={117535},
  year={(2025)},
  publisher={Elsevier}
}

@article{radermacher2014model,
  title={Model reduction in elastoplasticity: proper orthogonal decomposition combined with adaptive sub-structuring},
  author={Radermacher, Annika and Reese, Stefanie},
  journal={Computational Mechanics},
  volume={54},
  pages={677--687},
  year={(2014)},
  publisher={Springer}
}

@article{chaturantabut2010nonlinear,
  title={Nonlinear model reduction via discrete empirical interpolation},
  author={Chaturantabut, Saifon and Sorensen, Danny C},
  journal={SIAM Journal on Scientific Computing},
  volume={32},
  number={5},
  pages={2737--2764},
  year={(2010)},
  publisher={SIAM}
}

@article{kehls2023reduced,
  title={Reduced order modeling of structural problems with damage and plasticity},
  author={Kehls, Jannick and Kastian, Steffen and Brepols, Tim and Reese, Stefanie},
  journal={Proceedings in Applied Mathematics and Mechanics},
  volume={23},
  pages={e202300079},
  year={(2023)},
  publisher={Wiley Online Library}
}

@article{radermacher2016pod,
  title={POD-based model reduction with empirical interpolation applied to nonlinear elasticity},
  author={Radermacher, Annika and Reese, Stefanie},
  journal={International Journal for Numerical Methods in Engineering},
  volume={107},
  pages={477--495},
  year={(2016)},
  publisher={Wiley Online Library}
}

@article{farhat2014dimensional,
  title={Dimensional reduction of nonlinear finite element dynamic models with finite rotations and energy-based mesh sampling and weighting for computational efficiency},
  author={Farhat, Charbel and Avery, Philip and Chapman, Todd and Cortial, Julien},
  journal={International Journal for Numerical Methods in Engineering},
  volume={98},
  pages={625--662},
  year={(2014)},
  publisher={Wiley Online Library}
}

@article{an2008optimizing,
  title={Optimizing cubature for efficient integration of subspace deformations},
  author={An, Steven and Kim, Theodore and James, Doug},
  journal={ACM Transactions on Graphics},
  volume={27},
  pages={1--10},
  year={(2008)},
  publisher={ACM New York, NY, USA}
}

@article{bruno2013fracture,
  title={A fracture-ALE formulation to predict dynamic debonding in FRP strengthened concrete beams},
  author={Bruno, Domenico and Greco, Fabrizio and Lonetti, Paolo},
  journal={Composites Part B: Engineering},
  volume={46},
  pages={46--60},
  year={(2013)},
  publisher={Elsevier}
}

@article{ammendolea2025finite,
  title={Finite element modeling of dynamic crack branching using the moving mesh technique},
  author={Ammendolea, Domenico and Fabbrocino, Francesco and Leonetti, Lorenzo and Lonetti, Paolo and Pascuzzo, Arturo},
  journal={Engineering Fracture Mechanics},
  pages={111438},
  year={(2025)},
  publisher={Elsevier}
}

@article{ammendolea2025efficient,
  title={An efficient moving-mesh strategy for predicting crack propagation in unidirectional composites: Application to materials reinforced with aligned CNTs},
  author={Ammendolea, Domenico and Fabbrocino, Francesco and Leonetti, Lorenzo and Lonetti, Paolo and Pascuzzo, Arturo},
  journal={Composite Structures},
  volume={352},
  pages={118652},
  year={(2025)},
  publisher={Elsevier}
}

@article{ammendolea2023fatigue,
  title={Fatigue crack growth simulation using the moving mesh technique},
  author={Ammendolea, Domenico and Greco, Fabrizio and Leonetti, Lorenzo and Lonetti, Paolo and Pascuzzo, Arturo},
  journal={Fatigue \& Fracture of Engineering Materials \& Structures},
  volume={46},
  pages={4606--4627},
  year={(2023)},
  publisher={Wiley Online Library}
}

\end{document}